\documentclass[]{article}
\usepackage{fix-cm}
\usepackage{arxiv}
\usepackage[english]{babel}
\usepackage{url}
\usepackage{comment}
\usepackage{amsmath,amssymb}
\usepackage{mathtools}
\usepackage{graphicx}
\usepackage{caption}
\usepackage{floatrow}
\newfloatcommand{capbtabbox}{table}[][\FBwidth]
\usepackage{amsthm}
\usepackage{bm}
\usepackage{subfigure}
\usepackage{epstopdf}
\usepackage{float}
\usepackage[normalem]{ulem}
\usepackage{xcolor,colortbl}
\usepackage{enumitem}
\usepackage{todonotes}
\usepackage{hyperref}
\hypersetup{colorlinks=true, citecolor=blue}
\definecolor{abstract_background}{RGB}{225,245,225}%
\definecolor{Pink}{rgb}{1,0.8,0.9}
\definecolor{LightCyan}{rgb}{0.88,1,1}
\newcolumntype{g}{>{\columncolor{Pink}}c}

\newtheorem{definition}{Definition}[section]
\newtheorem{remark}[definition]{Remark}

\setenumerate[0]{label=\roman*)}
\setenumerate[1]{label=\arabic*.}
\setenumerate[2]{label=(\alph*)}

\newcommand{\shorttitle}{Solving Symmetry-Driven PDE dynamics with PINNs}
\title{Going with the Flow: \\
Solving for Symmetry-Driven PDE dynamics \\ with Physics-informed Neural Networks}

\author{
\textbf{
Michail E. Kavousanakis\textcolor{teal}{$^{1,*}$},
Gianluca Fabiani\textcolor{teal}{$^{2,3}$},
Anastasia Georgiou\textcolor{teal}{$^{3}$},
}
{}\\
\textbf{
Constantinos Siettos\textcolor{teal}{$^{4}$},
Panagiotis Kevrekidis\textcolor{teal}{$^{5}$},
Ioannis G.  Kevrekidis\textcolor{teal}{$^{3,6,7,}$\thanks{Corresponding authors, email: \texttt{mihkavus@chemeng.ntua.gr, yannisk@jhu.edu}}}}
}
\institute{
\textcolor{teal}{$^{(1)}$} School of Chemical Engineering, \emph{National Technical University of Athens}, 9 Iroon Polytechniou Street,\\ \hspace{0.41cm} Athens, 15772, Greece\\
\textcolor{teal}{$^{(2)}$} Hopkins Extreme Materials Institute, \emph{Johns Hopkins University}, Baltimore, 21218, MD, USA \hspace{1cm}\\
\textcolor{teal}{$^{(3)}$} Department of Chemical and Biomolecular Engineering, \emph{Johns Hopkins University}, Baltimore, 21218, MD, USA \\
\textcolor{teal}{$^{(4)}$}Dipartimento di Matematica e Applicazioni ‘‘Renato Caccioppoli", \emph{Universit\`a degli Studi di Napoli}\\ \hspace{0.41cm}\emph{Federico II}, Naples. 80126, Italy\\
\textcolor{teal}{$^{(5)}$} 
Department of Mathematics and Statistics and Department of Physics, \emph{University of Massachusetts}, \\
\hspace{0.41cm}Amherst, 01003, MA, USA\\
\textcolor{teal}{$^{(6)}$} Department of Applied Mathematics and Statistics, \emph{Johns Hopkins University}, Baltimore, 21218, MD, USA\\
\textcolor{teal}{$^{(7)}$} School of Medicine’s Dept. of Urology, \emph{Johns Hopkins University}, Baltimore, 21218, MD, USA
}

\begin{document}

\maketitle
\begin{abstract}
\colorbox{abstract_background}{\begin{minipage}{1\linewidth}
In the past, we have presented a systematic computational framework for analyzing self-similar and traveling wave dynamics in nonlinear partial differential equations (PDEs) by dynamically factoring out continuous symmetries such as translation and scaling. 
This is achieved through the use of time-dependent transformations—what can be viewed as dynamic pinning conditions—that render the symmetry-invariant 
solution stationary or slowly varying in rescaled coordinates. 
The transformation process yields a modified evolution 
equation coupled with algebraic constraints on the symmetry parameters, resulting in index-2 differential-algebraic equation (DAE) systems.
The framework accommodates both first-kind and second-kind self-similarity, and directly recovers the self-similarity exponents or wave speeds as part of the solution, upon
considering steady-state solutions in the rescaled
coordinate frame.
To solve the resulting high-index DAE systems, we employ Physics-Informed Neural Networks (PINNs), which naturally integrate PDE residuals and algebraic constraints into a unified loss function. 
This allows simultaneous inference of both the invariant solution and the transformation properties (such as the speed or the
scaling rate without the need for large computational domains, mesh adaptivity, or front tracking.
We demonstrate the effectiveness of the method on four canonical problems: (i) the Nagumo equation exhibiting traveling waves, (ii) the diffusion equation (1D and 2D) with first-kind self-similarity, (iii) the 2D axisymmetric porous medium equation showcasing second-kind self-similarity, and (iv) the Burgers equation, which involves both translational and scaling invariance.
The results demonstrate the capability of PINNs to effectively solve these complex PDE-DAE systems, providing a promising tool for studying nonlinear wave and scaling phenomena.
\end{minipage}}
\end{abstract}


\section{Introduction}
Many nonlinear partial differential equations (PDEs) arising in physics, biology, and engineering exhibit solutions characterized by symmetry-driven behavior, such as traveling waves, self-similar collapse, finite-time blow-up, rotating patterns, or radial expansions \cite{barenblatt1996scaling, whitham2011linear, merle2005blow, biktasheva2010computation, descalzi2013exploding, rowley2000reconstruction, aronson2001going, rowley2003reduction, siettos2003focusing, chen2004molecular, kavousanakis2007projective, kevrekidis2017infinity, wang2023asymptotic,LUSHNIKOV20101678,LUSHNIKOV20101678,Fibich2015,SulemSulem1999,KivsharPelinovsky2000}. These behaviors often reflect underlying continuous symmetries in the governing equations, such as translational, scaling, or rotational invariance. 
Exploiting these symmetries can simplify the analysis and reduce the complexity of the observed dynamics, but it typically requires reformulating the PDE by introducing additional terms 
-- such as convective or scaling components -- that account 
for the effects of translation, dilation, or rotation. These modifications effectively anchor the solution in a co-moving or rescaled frame, removing symmetry-induced drift or growth, but result in a transformed system that includes additional unknowns and algebraic constraints, often increasing the mathematical and numerical complexity of the problem \cite{aronson2001going}.
For example, traveling wave problems can be reformulated in a moving frame with unknown speed, while self-similar problems may involve dynamically rescaled coordinates, often governed by parameters, such as the rate of 
change of the width or of the amplitude, that must be inferred as part of the solution.

To uniquely determine the transformed solution, one often supplements the reformulation with template-based conditions or additional {\em ad hoc} constraints. These may include pointwise anchoring (e.g., fixing a value or derivative at a reference location such as an inflection point), orthogonality conditions, or minimization of integral functionals measuring the distance between the evolving solution and a prescribed template function. Such constraints help isolate the intrinsic dynamics from symmetry-related modes and guide the solution toward a 
unique identification
in the transformed frame. When applied systematically, these approaches typically yield systems of differential-algebraic equations (DAEs), where the algebraic part enforces normalization, symmetry breaking, or projection conditions tailored to the specific traveling or scaling behavior of interest \cite{aronson2001going}.

High-index DAEs pose significant challenges for traditional numerical solvers \cite{gear1971numerical, brenan1995numerical, gear1990introduction, shampine1999solving}. They are frequently stiff, sensitive to initial condition consistency, 
and require specialized integration techniques. These difficulties are compounded in the presence of multiple time scales, conservation laws, or quasi-steady-state approximations, where standard time-stepping methods may fail to capture the correct dynamics or preserve invariants. Additionally, when dealing with singularities or long-range traveling structures, classical PDE solvers require either large computational domains or adaptive mesh refinement -- both of which increase computational cost and can lead to numerical instabilities. These features are further
compounded in the realm of higher-dimensional problems.

Recent advances in scientific machine learning have opened new pathways for solving forward \cite{raissi2018numerical,han2018solving,samaniego2020energy,chan2019machine,fabiani2021numerical, fabiani2023parsimonious} and inverse problems \cite{fresca2020comprehensive,kovachki2021neural,lee2020coarse, lee2023learning, vlachas2022multiscale, fabiani2024task, fabiani2025randonets, fabiani2025enabling} involving differential equations.
Early works in the 1990s already explored the use of neural networks for approximating solutions \cite{lee1990neural,dissanayake1994neural,meade1994numerical, gerstberger1997feedforward, lagaris1998artificial} and identifying system dynamics \cite{krischer1993model, masri1993identification, chen1995universal, taprantzis1997fuzzy, gonzalez1998identification, siettos1999advanced}, 
setting the foundations for modern physics-informed approaches. Building on these ideas, Physics-Informed Neural Networks (PINNs) \cite{karniadakis2021physics, raissi2019physics,lu2021deepxde} have emerged as a promising framework, where the governing equations, along with initial, boundary, and possibly algebraic constraints, are embedded directly into the training of neural networks. This approach allows for the solution of both classical initial/boundary value problems and more complex scenarios involving data-driven discovery of unknown dynamics, with applications ranging from fluid mechanics to materials science \cite{de2021physics, cai2021physics, wu2018physics, goswami2022physics}.
Beyond solving classical boundary value problems, PINNs have been used in a wide range of scientific applications, including data-driven discovery of hidden dynamics \cite{raissi2018hidden}, inference of model parameters \cite{karniadakis2021physics}, and integration with sparse observations \cite{karniadakis2021physics, arzani2021uncovering}. The framework has also been extended to tasks such as feedback linearization \cite{alvarez2023discrete}, observer design \cite{alvarez2024nonlinear}, and identification of slow-invariant manifolds in multiscale systems \cite{patsatzis2024slow}.
Such efforts have been extended past 
continuum problems to nonlinear dynamical lattice
ones~\cite{SAQLAIN2023107498,shahab2025physicsinformedneuralnetworkshighdimensional}.

PINNs can be readily applied to differential-algebraic systems by incorporating both differential operators and algebraic constraints directly into the loss function. This approach avoids the need for index reduction procedures, specialized DAE solvers, or dedicated time-stepping schemes \cite{moya2023dae,lee2024robust}. Constraints arising from template conditions, normalization requirements, or even integral identities can be enforced seamlessly as
additional (soft) penalty terms in the training objective.
This flexibility can be particularly useful in capturing self-similar and solitonic behaviors, where traditional numerical solvers often face difficulties due to mesh adaptivity or domain deformation requirements. Recent applications include PINN-based modeling of spontaneous imbibition and self-similar flows \cite{deng2021application}, the identification 
of asymptotically self-similar solutions of
axisymmetric Euler flows~\cite{wang2023asymptotic},
the study of optical solitons with conservation constraints \cite{wu2022prediction}, and nonlinear wave equations with adaptive training strategies \cite{guo2023solving}. These works highlight the ability of PINNs to handle symmetry-reduced formulations in a unified and mesh-free setting.

In this work, we propose a PINN-based framework tailored to solving PDEs that exhibit traveling waves, self-similarity, and blow-up behavior through transformation into a DAE system. The methodology builds on the template-based representations formulated by a number of the 
present authors~\cite{rowley2003reduction,siettos2003focusing,Kavous2,CHAPMAN2022133396} 
to reduce the system while retaining the essential dynamics. We demonstrate that this approach avoids the need for large computational domains or adaptive meshing by working directly in the transformed coordinates, solving the full coupled DAE system using the PINN loss formulation. Warm-start training strategies are employed to handle the initial conditions and improve convergence in stiff regimes.

We validate our framework across several representative examples: (i) the Nagumo equation, which exhibits traveling fronts; (ii) the linear diffusion equation with self-similar decay; (iii) Burgers’ equation, combining both traveling and scaling features; and (iv) the porous media equation in one and two dimensions, exhibiting self-similarity and nonlinear radially-symmetric spreading. These cases illustrate the generality of the approach and highlight the potential of PINNs to serve as a robust tool for analyzing symmetry-driven PDE dynamics that naturally lead to high-index DAE formulations.

In the text, we will use $\partial_x$ to denote the partial derivative $\frac{\partial}{\partial x}$, for notational convenience.

\section{Problem Statement}
\label{sec:problem}
We consider the numerical solution of partial differential equations (PDEs) of the form
\begin{equation}
\label{eq:PDE1}
\partial_t u(\bm{x},t) = \mathcal{L}_{\bm{x}}(u),
\end{equation}
ideally posed on an unbounded spatial domain $\bm{x} \in \mathbb{R}^d$, approximated numerically over a sufficiently large finite domain.

Our interest is in PDEs whose solutions exhibit symmetry-driven structures, such as self-similar profiles and traveling waves. These patterns often reflect scaling invariance, translational invariance, or other continuous symmetries inherent to the PDE dynamics.

Often, a self-similar solution to an autonomous evolution PDE, as in Eq. \eqref{eq:PDE1}, in the variables $(\bm{x},t)\equiv(x_1, x_2, \dots, x_d, t)$ takes the form

\begin{equation}
|t - t^*|^\beta F\left( \frac{x_1}{|t - t^*|^{\alpha_1}}, \dots, \frac{x_d}{|t - t^*|^{\alpha_d}} \right),
\end{equation}

where the similarity exponents $\alpha_1, \dots, \alpha_d, \beta$, and the function $F$ are determined from the PDE and appropiate boundary and initial conditions. The time $t^*$ is either a known initial time (thus setting $t > t^*$), or an unknown critical time (e.g., blow-up time, thus setting $t < t^*$). When the system also respects time-translation invariance, we can always assume, without loss of
generality, $t^* = 0$.

To simplify the discussion, we begin by considering one-dimensional problems ($d=1$). The PDE operator $\mathcal{L}_x$ is assumed to be, in general, nonlinear and to satisfy a scaling property. Specifically, for all $A, B, C > 0$, there exist constants $a, b \in \mathbb{R}$ such that

\begin{equation}
\label{eq:PDE_scaling_property}
\mathcal{L}_{\bm{x}} \left( B f\left( \frac{x-c}{A} \right) \right) = A^a \, B^b \, \mathcal{L}_y(f(y)), \quad y = \frac{x-c}{A}.
\end{equation}

To capture symmetry-driven behaviors—such as co-traveling and self-similar dynamics—we consider a transformed ansatz $w(y,\tau)$ in rescaled space-time coordinates:

\begin{equation}
\label{eq:scaled_shifted_ansatz}
u(x,t) = B(\tau(t)) \, w\left( y(\tau)\equiv\frac{x - c(\tau(t))}{A(\tau(t))}, \tau(t) \right),
\end{equation}
where $A(\tau)$ and $B(\tau)$ are unknown scaling functions, $c(\tau)$ is a time-dependent translation, and $\tau(t)$ is an internal rescaled time variable, for which no a-priori
assumptions are made; the associated scaling is
expected to be inferred from the PDE dynamics.
The primary objective of this transformation is to reformulate the original PDE problem \eqref{eq:PDE1} into a form in which the dominant symmetry structures become stationary or slowly varying in time.
Substituting Eq. \eqref{eq:scaled_shifted_ansatz} into the PDE in Eq. \eqref{eq:PDE_scaling_property} yields the transformed equation

\begin{equation}
\label{eq:PDE_transformed}
    \biggl( \partial_{\tau} B \, w + B\, \partial_y w \, \frac{-\partial_{\tau} c-\partial_{\tau} A \, y}{A} + B \, \partial_{\tau} w\biggr) \, \partial_t \tau = A^a \, B^b \, \mathcal{L}_y(f(y)).
\end{equation}
Let us define the rescaled time $\tau(t)$ via
\begin{equation}
    \label{eq:rescaled_time_tau}
    \partial_t \tau = \sigma \, A^a \, B^{b-1},
\end{equation}
where $\sigma = \text{sgn}(t - t^*)$. Then, substituting $\partial_t \tau$ as in Eq. \eqref{eq:rescaled_time_tau} and dividing the transformed equation by B, we obtain

\begin{equation}
\label{eq:PDE_transformed_final}
\sigma \left( \frac{\partial_{\tau} B}{B} \, w -
\frac{\partial_{\tau} A}{A} \, y \, \partial_y w - \frac{\partial_{\tau}c}{A} \, \partial_y w + \partial_\tau w \right) = \mathcal{L}_y(w).
\end{equation}

This transformation does not by itself yield a closed system (it is underdetermined): the profile \( w(y, \tau) \) evolves according to the transformed PDE, while the scaling functions \( A(\tau) \), \( B(\tau) \), and \( c(\tau) \) remain unspecified. To close the system, additional algebraic conditions are required. These may arise from normalization conditions, conservation laws, or template-based algebraic relations.
These relations, combined with the transformed PDE, naturally lead to a system of differential-algebraic equations (DAEs).
In many cases, the number of effective free parameters can be reduced via further problem-specific manipulations. A common approach is to introduce the blow-up rates $G(\tau) = \partial_\tau A / A$ of the solution width and $C(\tau) = \partial_\tau B / B$ of the solution amplitude, which allow reformulating the problem in terms of these logarithmic growth rates.

Several strategies exist to define such algebraic closures, each reflecting different physical or geometric properties of the evolving solution. In the following, we briefly outline common approaches and their mathematical characteristics.


\section{Closure Strategies for the Transformed PDE System}
In this section, we review several approaches to close the transformed PDE system in Eq. \eqref{eq:PDE_transformed_final}

\paragraph{(i) Pinning Conditions.}
Pinning conditions fix geometric features of the solution profile $w(y,\tau)$ at a given spatial location $y_*$ in the rescaled coordinates. Typical examples include fixing the value $w(y_*) = w_*$, a zero-crossing, or the location of a maximum. These constraints serve to remove translation 
and/or scaling ambiguities by anchoring the solution in space. 

\paragraph{(ii) Conservation-Based Constraints.}
For PDEs with conserved quantities (e.g., mass, energy), these invariants can be enforced directly on $w$. Such constraints typically involve integrals of $w$, leading to integral algebraic relations, often of index-2.
For example, if the total mass is conserved,
\begin{equation}
\int_{-\infty}^{\infty} w(y, \tau) \, dy = K,
\end{equation}
where $K$ is a constant.


\paragraph{(iii) Template conditions.}
To close the system in Eq. \eqref{eq:PDE_transformed_final}, we can also introduce a template function $T$, typically chosen based on prior knowledge, e.g., the initial condition or an approximate steady-state shape. 
A common strategy involves minimizing the $L^2$-distance between the current profile $w(x)$ and a scaled-translated version of $T$, thus defining the mismatch loss, 

\begin{equation}
\label{eq:loss_template}
E(A, B, c) := \int_\Omega \left( w(x) - B \cdot T\left( \frac{x - c}{A} \right) \right)^2 dx,
\end{equation}

and impose the following conditions, that correspond to minimization of Eq. \eqref{eq:loss_template} 

\begin{equation}
\label{eq:loss_template_derivatives}
\frac{\partial E}{\partial A} \Big|_{A=1, B=1, c=0} = 0, \quad
\frac{\partial E}{\partial B} \Big|_{A=1, B=1, c=0} = 0, \quad
\frac{\partial E}{\partial c} \Big|_{A=1, B=1, c=0} = 0.
\end{equation}



Explicitly, the conditions in Eq. \eqref{eq:loss_template_derivatives} correspond to the following algebraic constraints:

\begin{equation}
\label{eq:template_conditions_explicit}
\int_\Omega (w - T) \, y \, \partial_y T \, dy = 0, \quad
\int_\Omega (w - T) \, T \, dy = 0, \quad
\int_\Omega (w - T) \, \partial_y T \, dy = 0.
\end{equation}

These provide three algebraic equations that constrain the instantaneous values of $A(\tau)$, $B(\tau)$, and $c(\tau)$, thus helping to close the transformed PDE system and yielding a DAE formulation. The extremization with
respect to $c$ provides information with respect to the
wave speed, while the other two determine the width
and amplitude information. In case these two are related
in their evolution, the relevant conditions may degenerate
into a single one for the ``blowup rate''.

\paragraph{Related works: method of slices, freezing method, and MN-dynamics.}
Several techniques have been introduced to eliminate continuous symmetries, such as translations and scalings, in the dynamics of PDEs. One classical approach is the \emph{method of slices}, introduced in \cite{rowley2003reduction}. This method defines a slice hyperplane that intersects the group orbit transversely, allowing the full system trajectory to be projected onto a reduced state space where the symmetry is factored out. The evolution in this reduced space is then coupled to reconstruction equations for the symmetry parameters. 
A related construction is the \emph{MN-dynamics framework} in \cite{aronson2001going}, where self-similar solutions are studied by introducing moving and rescaled coordinates. The authors derive equations for both the evolving profile and the parameters (e.g., speed, scale) by imposing additional constraints—often involving moment conditions or orthogonality, that regularize the dynamics. While similar in objective to the template-based method, MN-dynamics differs in the form of these constraints and in its treatment of the modulated profile. 
Such renormalization methods, already well-established
in the dispersive wave literature, e.g., through the
dynamic rescaling method and variants thereof~\cite{SulemSulem1999,Fibich2015}, have recently
been used to extract normal forms for the bifurcation
of collapsing solutions~\cite{Kavous2,Chapman_2024},
as well as to characterize the stability of the
collapsing solutions (as stationary states) in the
self-similar (so-called ``co-exploding'') frame.
The \emph{freezing method}, proposed in \cite{beyn2004freezing}, adopts a similar idea by transforming the PDE into a comoving and rescaled frame, yielding a differential-algebraic system. Phase conditions are imposed to "freeze" the solution in the transformed coordinates, typically by requiring orthogonality to the group tangent directions or by minimizing a suitable functional.

\paragraph{Index-reduction.}
One can arbitrarily select or combine any of the strategies described at points (i), (ii), or (iii) to impose a number of algebraic constraints matching the unknown transformation parameters (e.g., $A(\tau), B(\tau), c(\tau)$). All the algebraic constraints presented can close the transformed PDE system, in Eq \eqref{eq:PDE_transformed_final}, but result in an index-2 DAE, since they involve only linear functionals of $w$ (either at specific spatial points or over domains) and do not explicitly depend on the transformation parameters. Consequently, two derivatives with respect to $\tau$ are needed to obtain an explicit ODE system.
A common index-reduction technique involves differentiating the constraint once with respect to $\tau$, reducing the system to index-1. Let the constraint take the generic form $L_{cons}(w) + G_{cons} = 0$, where $L_{cons}$ is a linear operator and $G_{cons}$ is independent of $w$ and $\tau$. Differentiating yields

\begin{equation}
\partial_\tau \left( L_{cons}(w) + G_{cons} \right) = L_{cons}(\partial_\tau w) = 0.
\end{equation}

Then, substituting the transformed PDE expression for $\partial_\tau w$, as in Eq. \eqref{eq:PDE_transformed_final},

\begin{equation}
L_{cons}\big(\partial_\tau w\big) = L_{cons}\bigg( \sigma \mathcal{L}_y(w)
     - \frac{\partial_\tau B}{B} \, w
     - \frac{\partial_\tau A}{A} \, y \, \partial_y w
     - \frac{\partial_\tau c}{A} \, \partial_y w \bigg)=0,
\end{equation}

the differentiated constraint $L_{cons}(\partial_\tau w) = 0$ becomes a relation involving explictly the unknowns $A, B, c$ and their derivatives. This yields an index-1 DAE system that is more amenable to numerical integration.

In the context of Physics-Informed Neural Networks, however, this index-reduction is not strictly necessary. The original index-2 formulation can be enforced directly through the loss function, which simultaneously penalizes both the PDE residual and the algebraic constraints. This flexibility allows the method to bypass the derivation of auxiliary evolution equations, which are typically required by conventional solvers.

\section{Numerical solution of DAEs with PINNs}
In this section and for the completeness of presentation, we first give a very brief introduction to the basic concept of physics-informed machine learning for the solution to differential equations in the form of PDEs/DAEs.

To approximate solutions of time-dependent PDEs involving evolving symmetry parameters, we consider a formulation combining differential and algebraic components. Let the main PDE be

\begin{equation}
\partial_{\tau} u(x,t) = \mathcal{L}_x\left(u; \bm{p}(t)\right), \quad x \in \Omega \subset \mathbb{R}^d,
\end{equation}

with boundary conditions

\begin{equation}
\mathcal{B} u(x,t) = g(x,t), \quad x \in \partial \Omega.
\end{equation}

Here, $\bm{p}(t)\in\mathbb{R}^m$ denotes a set of time-dependent parameters (e.g., scaling or translation functions such as $A(t), B(t), c(t)$), which appear explicitly in the transformed equation or influence the dynamics implicitly via the chosen frame.

To close the system, we assume that $\bm{p}(t)$ must satisfy additional algebraic constraints:

\begin{equation}
f_j(u(\cdot,t), \bm{p}(t)) = 0, \quad j=1,\dots,m,
\end{equation}

where each $f_j$ may encode conservation laws, moment conditions, or template-matching criteria.

We introduce two neural networks: $\mathcal{N}_u(x,t)$ to approximate the solution $u(x,t)$, and $\mathcal{N}_{\bm{p}}(t)$ to approximate the evolving parameters $\bm{p}(t)$.
Let $\{x_i\}_{i=1}^{n_x} \subset \Omega$, $\{x_{\ell}^{\partial}\}_{\ell=1}^{n_{\partial}} \subset \partial \Omega$, and $\{t_k\}_{k=1}^{n_t}$ be sampling points in space and time.
The residual of the differential operator is incorporated into a loss function:

\begin{equation}
\mathcal{E}_\text{pde} := \sum_{i=1}^{n_x} \sum_{k=1}^{n_t} \left| \partial_t \mathcal{N}_u(x_i,t_k) - \mathcal{L}_x\left( \mathcal{N}_u; \mathcal{N}_p(t_k) \right) \right|^2,
\end{equation}

along with algebraic residuals

\begin{equation}
\mathcal{E}_\text{alg} := \sum_{j=1}^{m} \sum_{k=1}^{n_t} \left| f_j\left(\mathcal{N}_u(\cdot,t_k), \mathcal{N}_p(t_k) \right) \right|^2,
\end{equation}

boundary residuals

\begin{equation}
\mathcal{E}_\text{bc} := \sum_{\ell=1}^{n_{\partial}} \sum_{k=1}^{n_t} \left| \mathcal{B} \mathcal{N}_u(x_{\ell}^{\partial}, t_k) - g(x_{\ell}^{\partial}, t_k) \right|^2.
\end{equation}

and initial condition residuals:

\begin{equation}
\mathcal{E}_\text{ic} := \sum_{i=1}^{n_{x}} \left| u(x_{i}, t_0) - u_0(x_{i}) \right|^2.
\end{equation}

The total loss is given by

\begin{equation} \label{eq:lossterm}
\mathcal{E}_\text{total} = \lambda_\text{pde} \mathcal{E}_\text{pde} + \lambda_\text{alg} \mathcal{E}_\text{alg} + \lambda_\text{bc} \mathcal{E}_\text{bc} + \lambda_\text{ic} \mathcal{E}_\text{ic}
\end{equation}

with weights $\lambda_\text{pde}, \lambda_\text{alg}, \lambda_\text{bc}, \lambda_\text{ic}$ to balance contributions. Backpropagation/ Automatic differentiation can provide the required time and spatial derivatives, as well as the gradients w.r.t trainable parameters of the PINN. This leads to a physics-informed neural approximation of a system of DAEs, where both the field variable and its symmetry parameters evolve consistently with the governing laws and closure conditions.

Ultimately, in the context of symmetry-reduced formulations—such as co-moving frames for traveling waves or self-similar scaling frames—additional algebraic constraints may arise from the coordinate transformation. These include template conditions, normalization constraints, or integral identities used to fix the degrees of freedom introduced by symmetry. PINNs can naturally incorporate these terms as soft penalties into the loss function, enabling the consistent solution of mixed PDE-DAE systems arising in the evolution of rescaled or traveling profiles.

\section{Numerical Results}
In this section, we showcase the application of the PINN-assisted MN framework across four distinct case studies: 
one involving a traveling wave solution (Nagumo equation), two exhibiting self-similar behavior (1D diffusion and 2D axisymmetric porous medium equation), and one combining both scaling and translational invariance (Burgers equation).
%
For each case, we formulate an index-2 DAE system and demonstrate how the PINN-assisted methodology can be used to numerically recover both the invariant (translational or self-similar) profile and the corresponding wave speed or self-similar exponents within a dynamically rescaled coordinate system.
All simulations were performed in \texttt{MATLAB} on a laptop with an $11^{th}$ Gen Intel(R) Core(TM) i7-111800H @2.30CPU GHz processor, 32 GB RAM, and Windows 11 Home operating system. The neural networks were implemented using \texttt{dlnetwork} and \texttt{dlarray} objects to enable automatic differentiation and GPU acceleration when available.
The architecture of the neural networks follows a simple feedforward structure with hyperbolic tangent activations. Each model is defined using the functions \texttt{featureInputLayer}, \texttt{fullyConnectedLayer}, \texttt{tanhLayer}.
Training was performed using a quasi-Newton method (L-BFGS) via a custom \texttt{lbfgsupdate} routine, suitable for low-batch, physics-informed settings.

\begin{remark}
    As a practical strategy to improve convergence, we introduced a warm-up phase prior to full training. Specifically, the network is first trained to approximate the initial condition alone, extended trivially over all times (i.e., enforcing that the initial condition remains constant across time). This preconditioning step helps the network establish a suitable parameter initialization before minimizing the full physics-informed loss, leading to more stable convergence in the subsequent optimization.
\end{remark}


\subsection{Case study 1: Traveling wave in the Nagumo PDE}

We first illustrate the proposed methodology using a translationally invariant problem: the Nagumo equation, given by
\begin{equation}\label{eq:nagumo}
    \partial_t u = \partial^2_{xx} u + u(1-u)(u-a).
\end{equation}
\begin{figure}[ht!]
    \centering
    \subfigure[]{\includegraphics[width=0.45\linewidth]{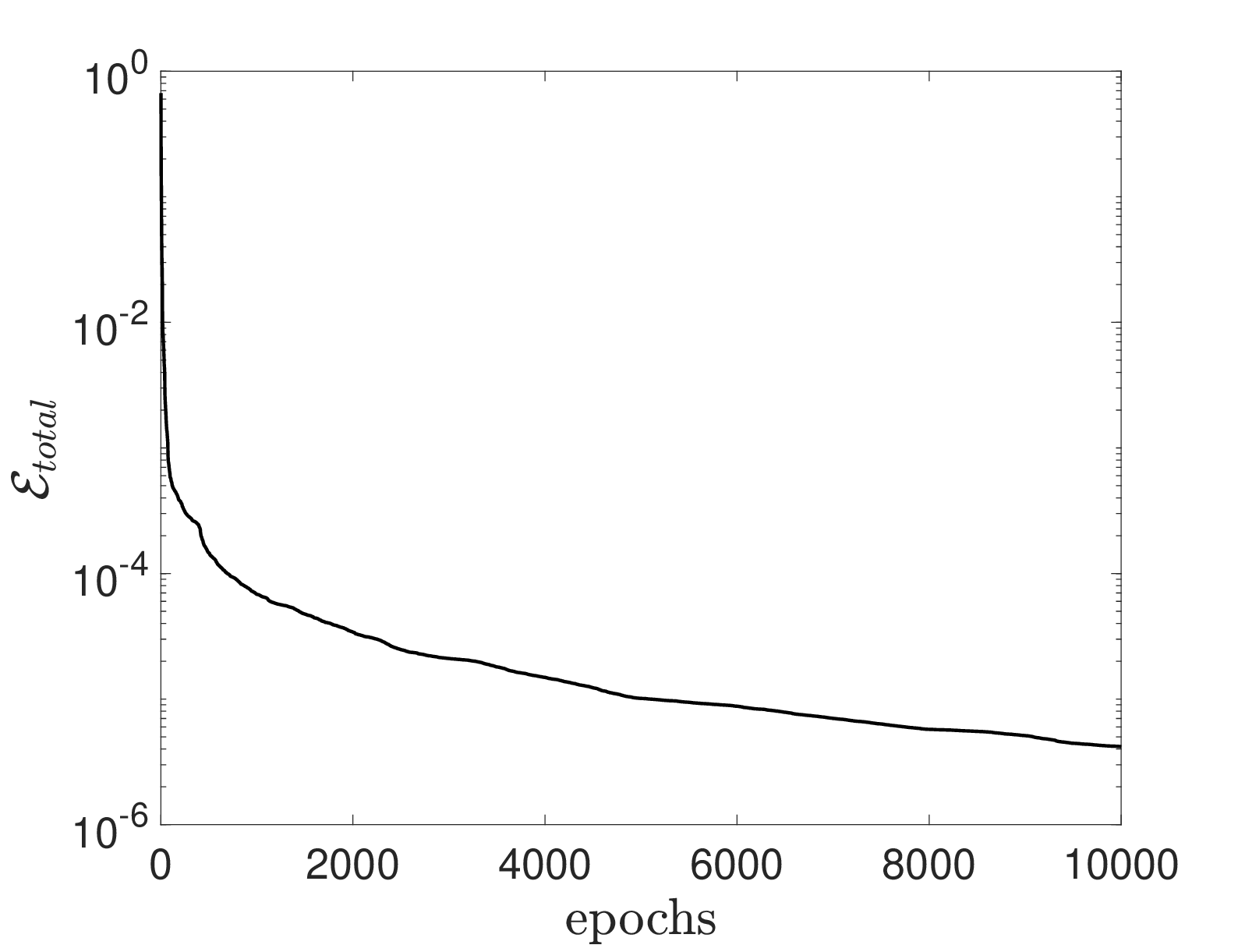} \label{fig:nagumo_case_a}}
    \subfigure[]{ \includegraphics[width=0.45\linewidth]{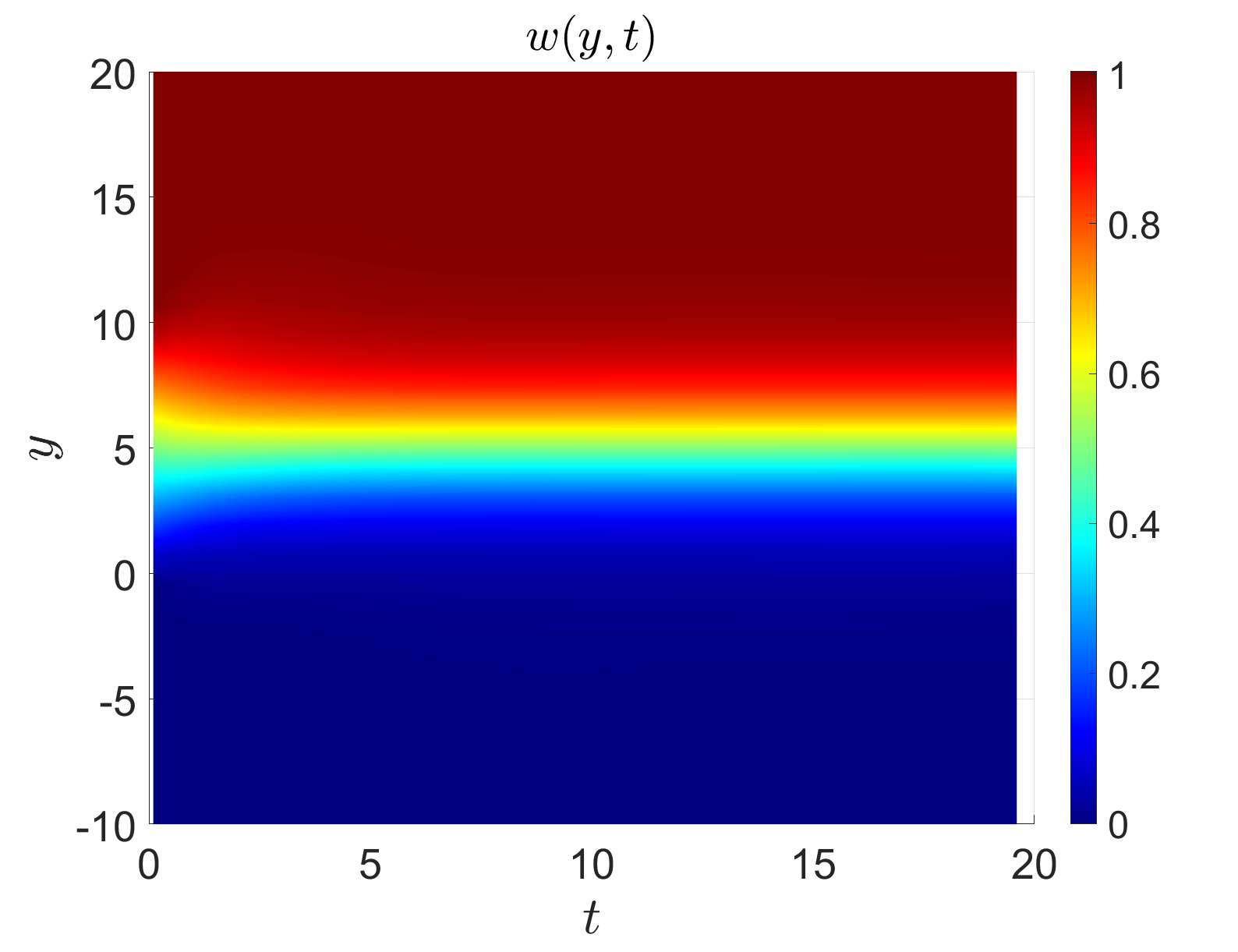} \label{fig:nagumo_case_b}}
    \hfill
    \subfigure[]{\includegraphics[width=0.45\linewidth]{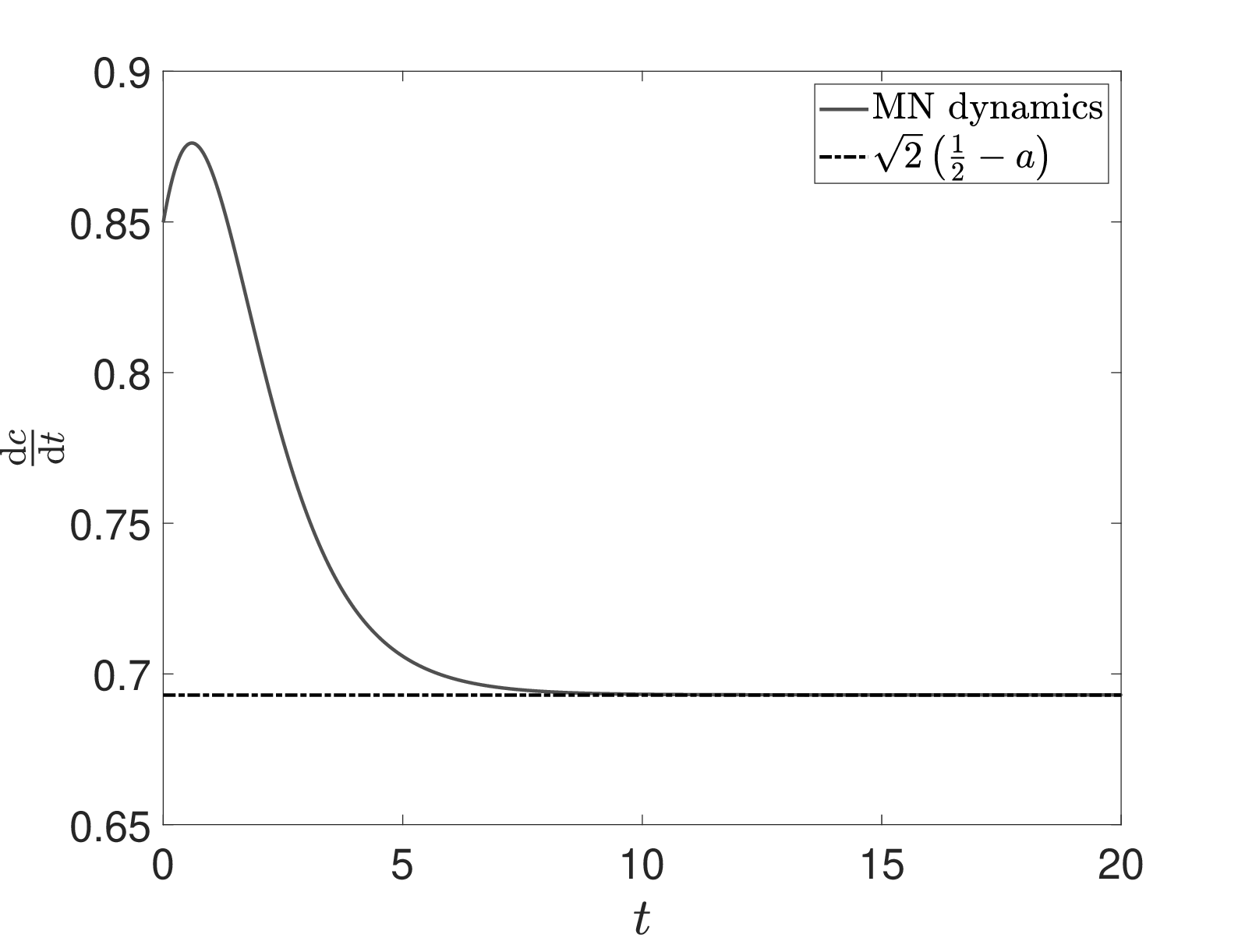} \label{fig:nagumo_case_c}}
    \subfigure[]{\includegraphics[width=0.45\linewidth]{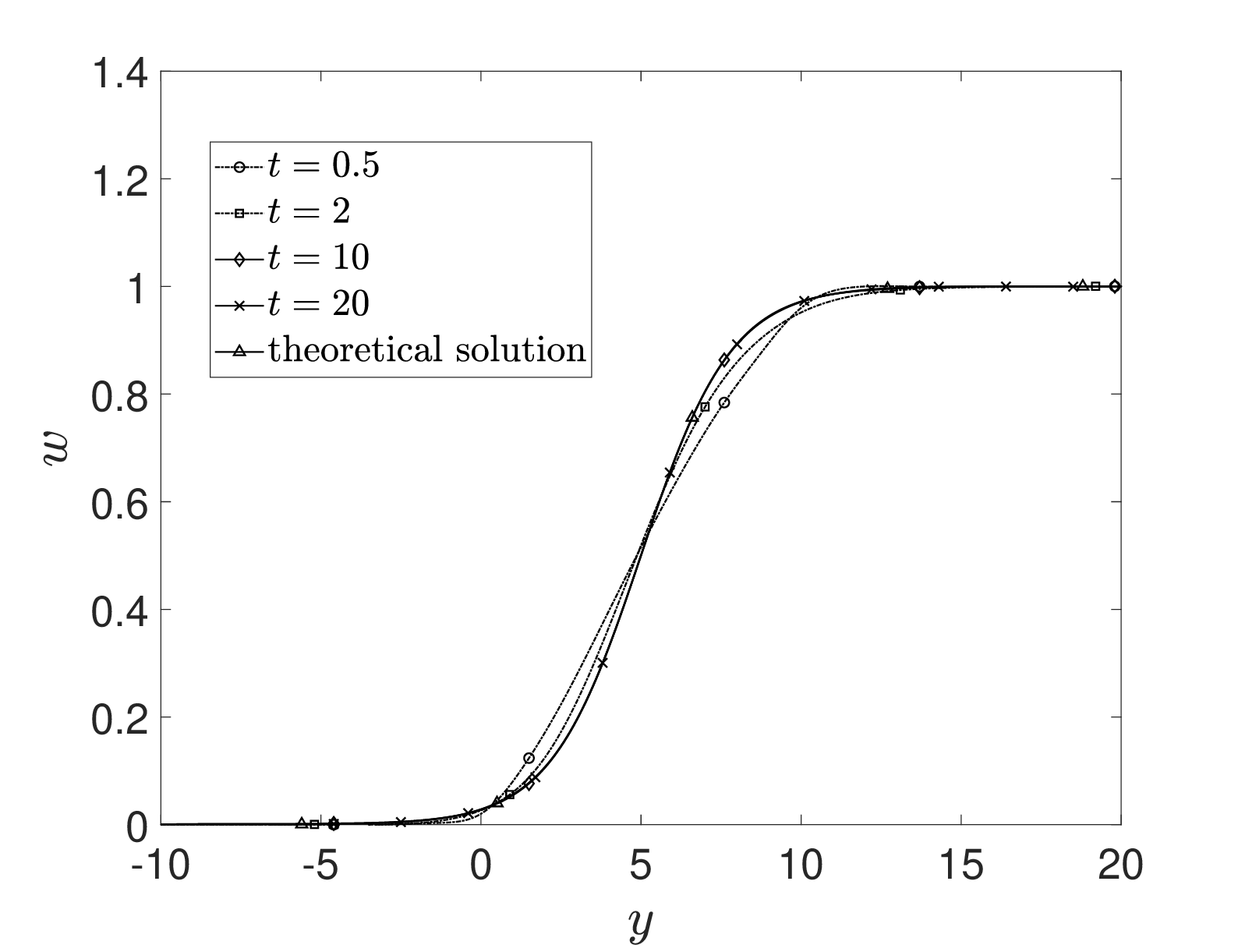} \label{fig:nagumo_case_d}}
    \caption{
    PINN prediction of the MN dynamics for the co-moving Nagumo equation in Eq \eqref{eq:nagumo}.
    (a) Convergence of the total loss $\mathcal{E}_{\text{total}}$ during training.
    (b) Evolution of the rescaled solution $w(y,t)$ in co-moving coordinates.
    (c) Predicted instantaneous velocity $\frac{\mathrm{d}c}{\mathrm{d}t}$ for $a = 0.01$, approaching the theoretical value (black dashed-dotted line).
    (d) Snapshots of the rescaled (co-moving) solution $w$ at various time instances.
    At long times ($t > 10$), the PINN-predicted solution $w$ converges to a stationary profile that shows excellent agreement with the theoretical expression in Eq.~(\ref{eq:nagumotheory}), up to a spatial shift of approximately 5 units in $x$.
    \label{fig:nagumo_case}}
\end{figure}
\noindent The Nagumo equation is a classical example of a parabolic PDE that supports traveling wave solutions.
An explicit traveling wave solution takes the form $u(x,t)=w(x-c)$, where:
\begin{equation} \label{eq:nagumotheory}
    w(x) = \left[ 1+ \textrm{exp}\left( -\frac{x}{\sqrt{2}} \right) \right]^{-1}, \quad \frac{\textrm{d}c}{\textrm{d}t}=-\sqrt{2} \left( \frac{1}{2}-a \right).
\end{equation}

\noindent To capture the traveling wave behavior, we dynamically construct a co-moving (co-traveling) frame using a template condition to recover the instantaneous wave speed.
In this framework, the solution is written as: 
\begin{equation} \label{eq:solcomoving}
    u(x,t)=w(x-c(t),t).
\end{equation}

\noindent Substituting Eq.~(\ref{eq:solcomoving}) into Eq.~(\ref{eq:nagumo}) and using the translational invariance of the operator: $\mathcal{D}_x (u) \equiv  \partial^2_{xx} u + u(1-u)(u-a) $, we obtain the equation in the co-moving coordinate, $y \equiv x-c(t)$:

\begin{equation} \label{eq:nagumo_MN}
    \partial_t w = \partial^2_{yy} w + w(1-w)(w-a) - \frac{\textrm{d}c}{\textrm{d}t} \partial_yw.
\end{equation}

The additional degree of freedom $\frac{\textrm{d}c}{\textrm{d}t}$ is determined using an algebraic-template condition.
Specifically, we enforce that the solution in the co-moving frame remains as close as possible to a predefined template function $T(y)$ by imposing the constraint: 

\begin{equation} \label{eq:tempnagumo}
    \int_{y_{min}}^{y_{max}} \left( w(y,t) - T \right) \partial_y T \textrm{d}y = 0. 
\end{equation}

For this case-study, Eq.~(\ref{eq:nagumo_MN}) is solved in the domain $y \in [-30,30]$ with Neumann boundary conditions and parameter $a=0.01$.
The initial condition is given by:

\begin{equation} \label{eq:nagumoIC}
    w(y,0) = \left\{ 
    \begin{array}{ll}
         0, & y<0  \\
         y/10, & 0 \leq y \leq 10 \\
         1, & y>10
    \end{array}
    \right.
\end{equation}

We also use this initial condition as the template function, i.e., $T(y) = w(y,0)$.

Two neural networks are employed: 
\begin{itemize}
    \item $\mathcal{N}_w(y,t)$: Takes $y$ and $t$ as inputs and approximates $w(y,t)$. 
    It has two hidden layers with 20 neurons each. 
    \item $\mathcal{N}_p(t)$: Takes $t$ as input and approximates the instantaneous wave speed $\frac{\textrm{d}c}{\textrm{d}t}$. 
    It consists of one hidden layer with 5 neurons. 
\end{itemize}

Both networks use the hyperbolic tangent as the activation function. 
They are trained simultaneously by minimizing a total loss function following the general form of Eq.~(\ref{eq:lossterm}), with: $\mathcal{E}_\text{alg} \equiv \int_{y_{min}}^{y_{max}} \left( w(y,t) - T \right) \partial_y T \textrm{d}y $. 
Training is performed using the L-BFGS optimizer, which ensures efficient convergence of the total loss $\mathcal{E}_{total}$ (see Figure~\ref{fig:nagumo_case_a}).
For $\mathcal{N}_w(y,t)$, we use 119,800 collocation in the $y-t$ domain  (599 equidistributed points in $y \in (-30,30)$ and 200 equally distributed points in $t \in (0,20]$.
For $\mathcal{N}_p(t)$, we use 201 equally spaced points in $t \in [0,20]$.
The spatio-temporal solution is depicted in Figure~\ref{fig:nagumo_case_b}.
Figure~\ref{fig:nagumo_case_c} displays the predicted instantaneous velocity $\frac{\textrm{d}c}{\textrm{d}t}$, which converges to approximately $0.69305$ at long $t$.
This closely matches the theoretical value $\sqrt{2} \left( \frac{1}{2}-0.01 \right) = 0.692965$.
As shown in Figure~\ref{fig:nagumo_case_d}, the template condition enforces the solution $w(y,t)$ in the co-moving frame to become stationary in shape, converging to a traveling wave form.
The comparison between $w$ at long $t$ ($t=20$) with the theoretical profile (Eq.~(\ref{eq:nagumotheory})) in Figure~\ref{fig:nagumo_case_d} also shows  an almost excellent agreement.

\subsection{Case study 2: Self-similarity in the Diffusion PDE}

\subsubsection{1D Diffusion}
As a second illustrative example, we consider the one-dimensional diffusion equation: 
\begin{equation} \label{eq:diffusion}
    \partial_t u = \partial^2_{xx} u,
\end{equation}

\begin{figure}[ht!]
    \centering
    \subfigure[]{
    \includegraphics[width=0.45\linewidth]{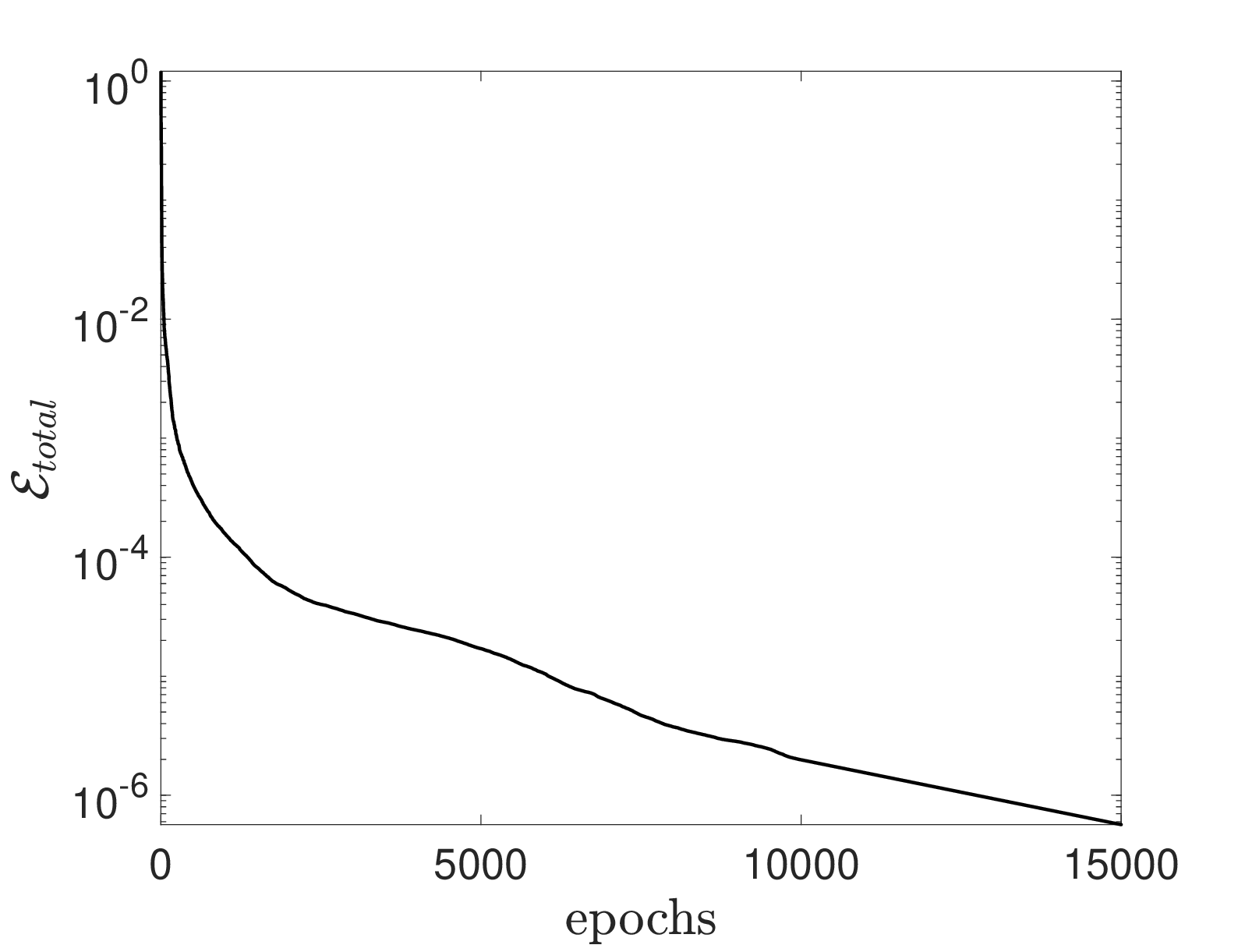}\label{fig:diffusionMN_a}}
    \subfigure[]{
    \includegraphics[width=0.45\linewidth]{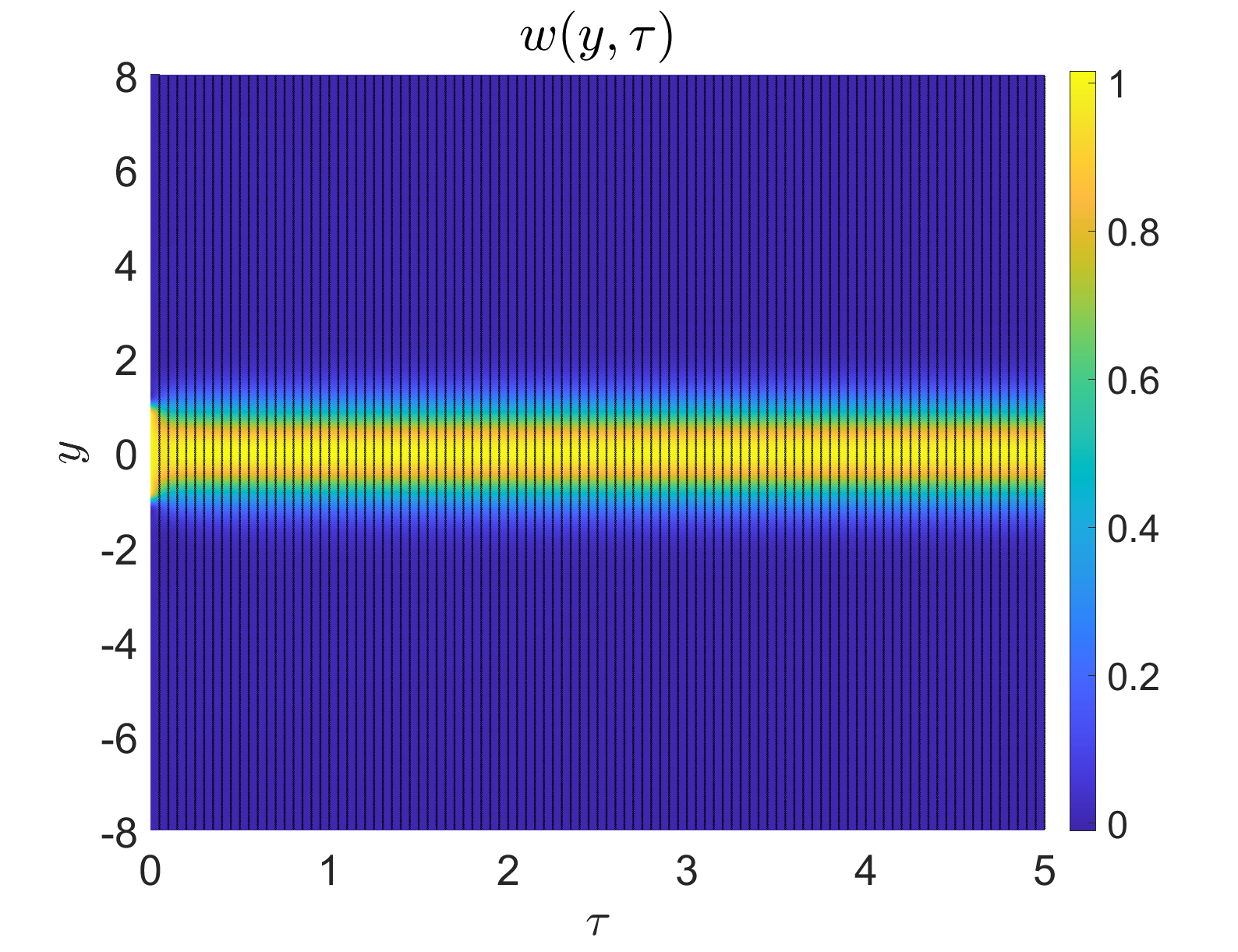}
    \label{fig:diffusionMN_b}
    }
    \subfigure[]{\includegraphics[width=0.45\linewidth]{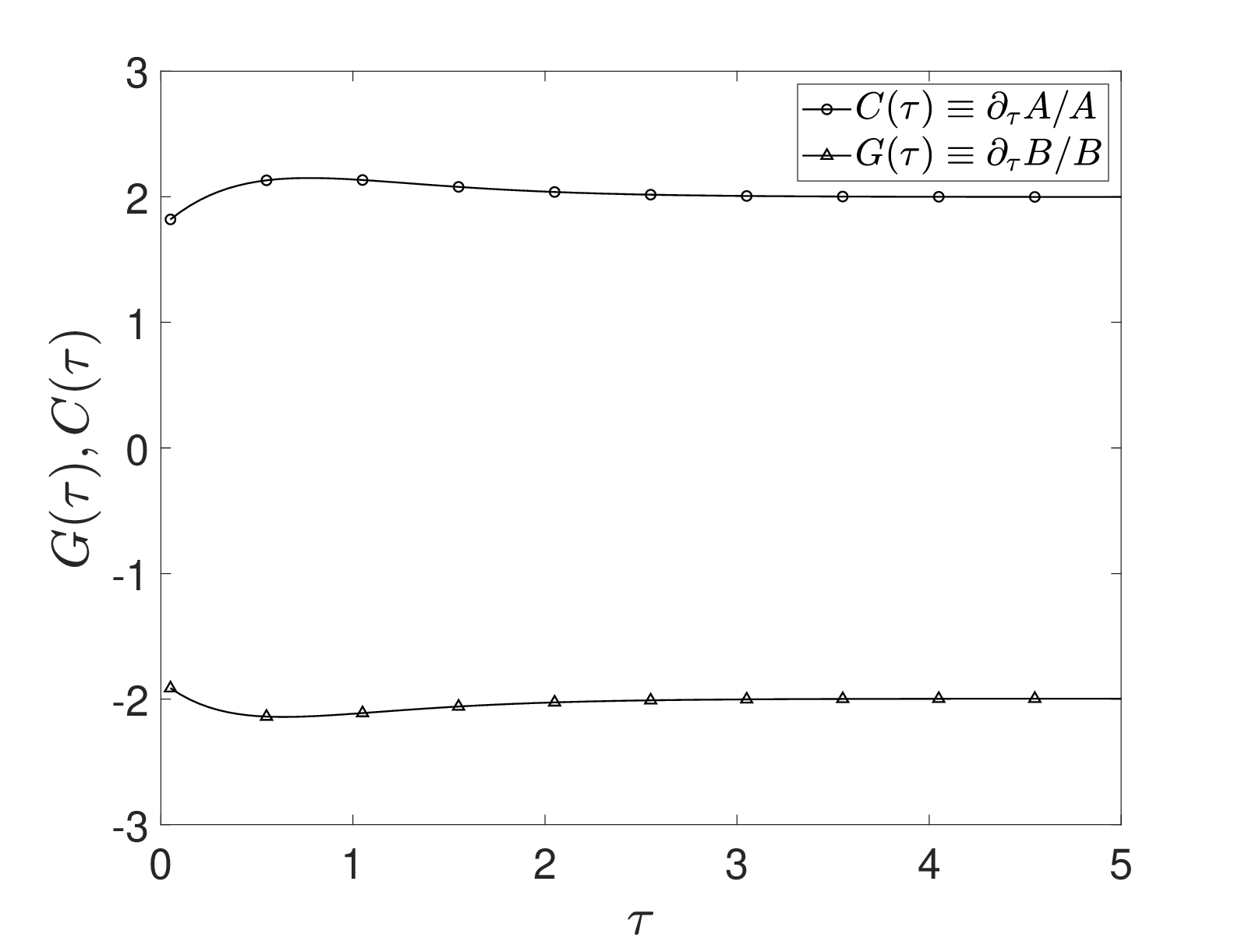}
    \label{fig:diffusionMN_c}
    }
    \subfigure[]{\includegraphics[width=0.45\linewidth]{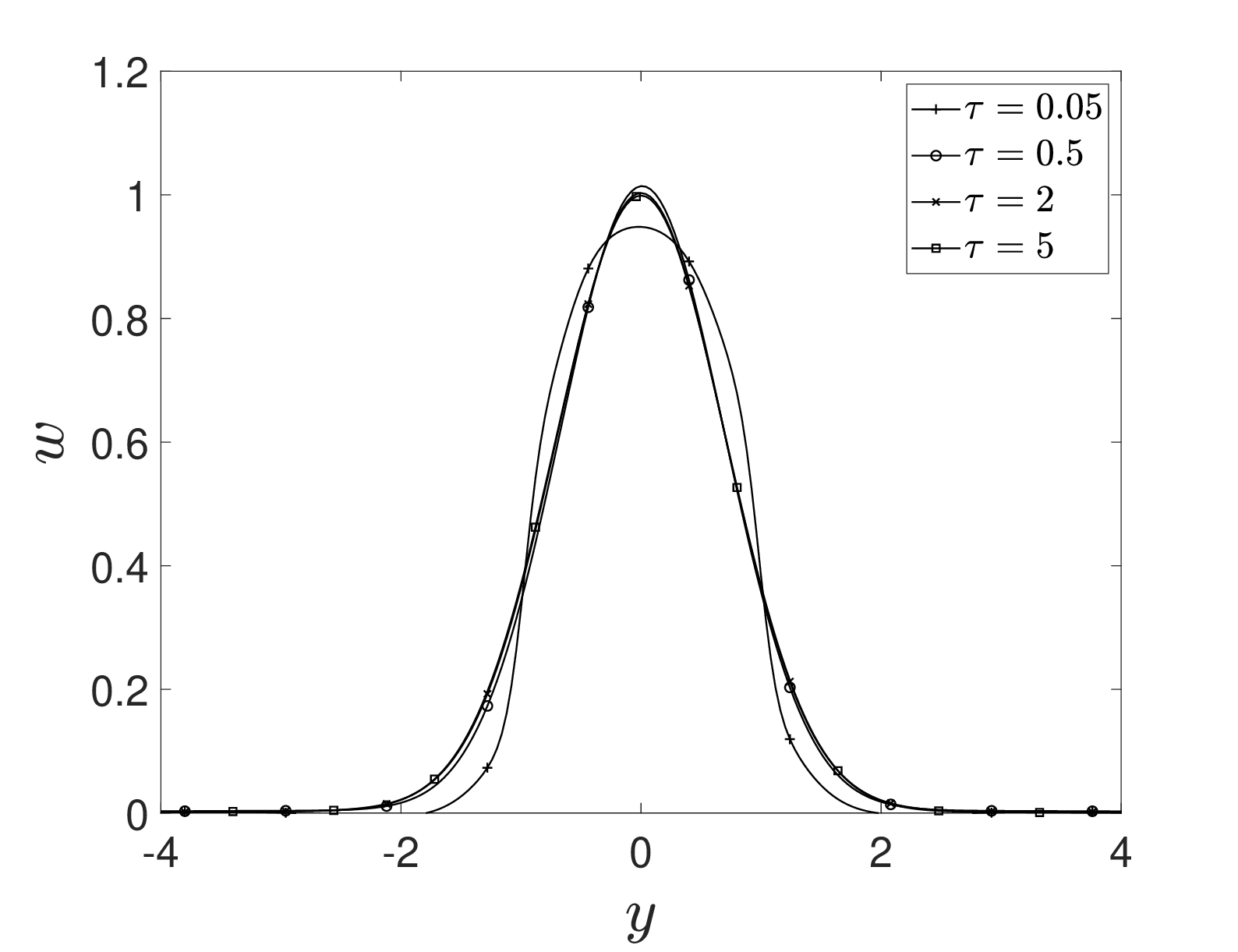}
    \label{fig:diffusionMN_d}
    }
    \caption{PINN prediction of the MN dynamics for the 1D diffusion equation in Eq \eqref{eq:diffusion}.
    (a) Convergence of the total loss $\mathcal{E}_{loss}$ during training for the  1D diffusion equation 
    in the transformed frame.
    (b) Spatio-temporal evolution of the rescaled solution $w(y,\tau)$. 
    The neural network $\mathcal{N}_w(y,\tau)$ is trained using 39,900 collocation points in the $y-\tau$ domain (indicated by black dots).
    (c) Evolution of scaling rates $G$ and $C$ over rescaled time, $\tau$.
    As $\tau \to \infty$, the ratio $C/G$ converges to a steady-state value, $\lim_{\tau \to \infty} {G/C} = -1$, which is used to derive the self-similar exponents of the diffusion equation.
    (d) Snapshots of the rescaled solution, $w(y,\tau)$ at different $\tau$ instances.
    Convergence to the final self-similar profile is effectively achieved by $\tau \approx 2$. 
    \label{fig:diffusionMN}}
\end{figure}

\noindent the right-hand side operator $\mathcal{D}_x (u) \equiv \partial^2_{xx} u $ satisfies the following scaling relation: 
\begin{equation} \label{eq:scalingdiffusion}
    \mathcal{D}_x \left ( B w \left(\frac{x}{A} \right) \right) = BA^{-2}\mathcal{D}_y w, \quad \textrm{ with } y=\frac{x}{A}
\end{equation}

In our dynamic renormalization framework, discussed
earlier, e.g., in~\cite{Kevrekidis2017}, we introduce the ansatz:
\begin{equation}
\label{eq:diffusion_ansatz}
u(x,t)=B(\tau) w \bigg(\frac{x}{A(\tau)}, \tau(t)\bigg),
\end{equation}

\noindent which, upon substitution of Eq.~(\ref{eq:diffusion_ansatz}) into Eq.~(\ref{eq:diffusion}) and using Eq.~(\ref{eq:scalingdiffusion}) yields: 

\begin{equation}\label{eq:diffusionMN}
    \partial_\tau w + \frac{\partial_\tau B}{B} \, w -\frac{\partial_\tau A}{A} \, y \, \partial_y w = \partial^2_{yy}w, \quad \textrm{ with } y=\frac{x}{A(\tau)}, \textrm{ and }  \partial_t \tau = A^{-2}. 
\end{equation}

\noindent To determine the scaling rates, 
$G(\tau) \equiv \frac{\partial_\tau B}{B}$ and $C(\tau) \equiv \frac{\partial_\tau A}{A}$, we impose two additional algebraic template conditions.
These conditions arise from minimizing the functional: 
\begin{equation} \label{eq:origtemplate}
    E \equiv \int_{y_{min}}^{y_{max}} \left( w - B T (\frac{\zeta}{A})^2 \right) \textrm{d}\zeta,
\end{equation}

\noindent when $A=B=1$.
This leads to the following constraints:
\begin{equation} \label{eq:template_diffusion}
    \int_{y_{min}}^{y_{max}} \left( w-T \right) T \textrm{d}y = 0 \quad \textrm{ and } \quad \int_{y_{min}}^{y_{max}} \left( w-T \right) \, y \, \partial_yT \, \textrm{d}y = 0,
\end{equation}

\noindent where $T$ is a template function. 
Here, $T(y)=\exp(-y^2)$ is chosen as the template function.
The initial condition is a smoothed pair of step functions
\begin{equation} \label{eq:initialdiffusion}
    w(y,0)= \frac{1}{2} \left( \tanh{ \frac{y+1}{0.2}} - \tanh{ \frac{y-1}{0.2}} \right).
\end{equation}

\noindent The spatial domain is $y \in [-8,8] $ with zero-flux (Neumann) boundary conditions at both ends. 

The neural network $\mathcal{N}_w(y,\tau)$, approximating $w(y,\tau)$, consists of two hidden layers with 40 neurons each and uses hyperbolic tangent activation functions. 
It admits as inputs the rescaled spatial and temporal variables $y$ and $\tau$.
The auxiliary network $\mathcal{N}_p(\tau)$, tasked with predicting the blow-up rates $G(\tau)$, and $C(\tau)$, takes as input only the rescaled time, $\tau$, and has a single hidden layer with 5 neurons.
The overall loss-function minimized during training includes: $\mathcal{E}_{pde}$ from Eq.~(\ref{eq:diffusionMN}), $\mathcal{E}_{bc}$ from boundary conditions, $\mathcal{E}_{ic}$ from the initial condition Eq.~(\ref{eq:initialdiffusion}), and $\mathcal{E}_{alg}$ from the algebraic template constraints in Eq.~(\ref{eq:template_diffusion}).
Training is performed using the L-BFGS optimizer and the convergence of the total loss $\mathcal{E}_{total}$ is shown in Figure~\ref{fig:diffusionMN_a}.
For $\mathcal{N}_w(y,t)$, we use 39,900 collocation points in the $y-t$ domain.  (399 equidistributed points in $y \in (-8,8)$ and 100 equally distributed points in $t \in (0,5]$.
For $\mathcal{N}_p(t)$, we use 100 equally spaced points in $t \in [0,5]$.
The spatio-temporal solution is depicted in Figure~\ref{fig:diffusionMN_b}.
The blow-up rates $G(\tau)$ and $C(\tau)$ are illustrated in
Figure~\ref{fig:diffusionMN_c}.
The rescaled solution converges to a stationary self-similar profile after $\tau \approx 2$, while the blow-up rates stabilize at $G_{steady}=-1.99699$ and $C_{steady}=1.99881$.
The ratio of these rates allows estimation of the self-similar exponents of the diffusion equation. 

Indeed, this expression leads to the corresponding
evolution of $B=B_0 e^{G_{steady} \tau}$
and $A=A_0 e^{C_{steady} \tau}$,
while this, in turn through the last of Equations~(\ref{eq:diffusionMN}), leads to 
$e^{C_{steady} \tau}=\frac{\sqrt{2 C_{steady}}}{A_0}
 (t+t_0)^{1/2}$, where $t_0$ is a non-negative
 integration constant. Using this derived transformation
 between $\tau$ and the original time $t$
 in the equations for
 $B$ and $A$, we infer that
 $A \propto (t+t_0)^{1/2}$, while 
 $B \propto (t+t_0)^{-G_{steady}/(2 C_{steady})}$
 which is $\approx -1/2$, per the computational
 results above.

\subsubsection{2D diffusion PDE}
For the two-dimensional diffusion equation: 

\begin{equation} \label{eq:diffusion2d}
    \partial_tu = \partial_{xx}^2 u + \partial_{yy}^2u, 
\end{equation}

\begin{figure}[ht!]
    \centering
    \subfigure[]{\includegraphics[width=0.45\linewidth]{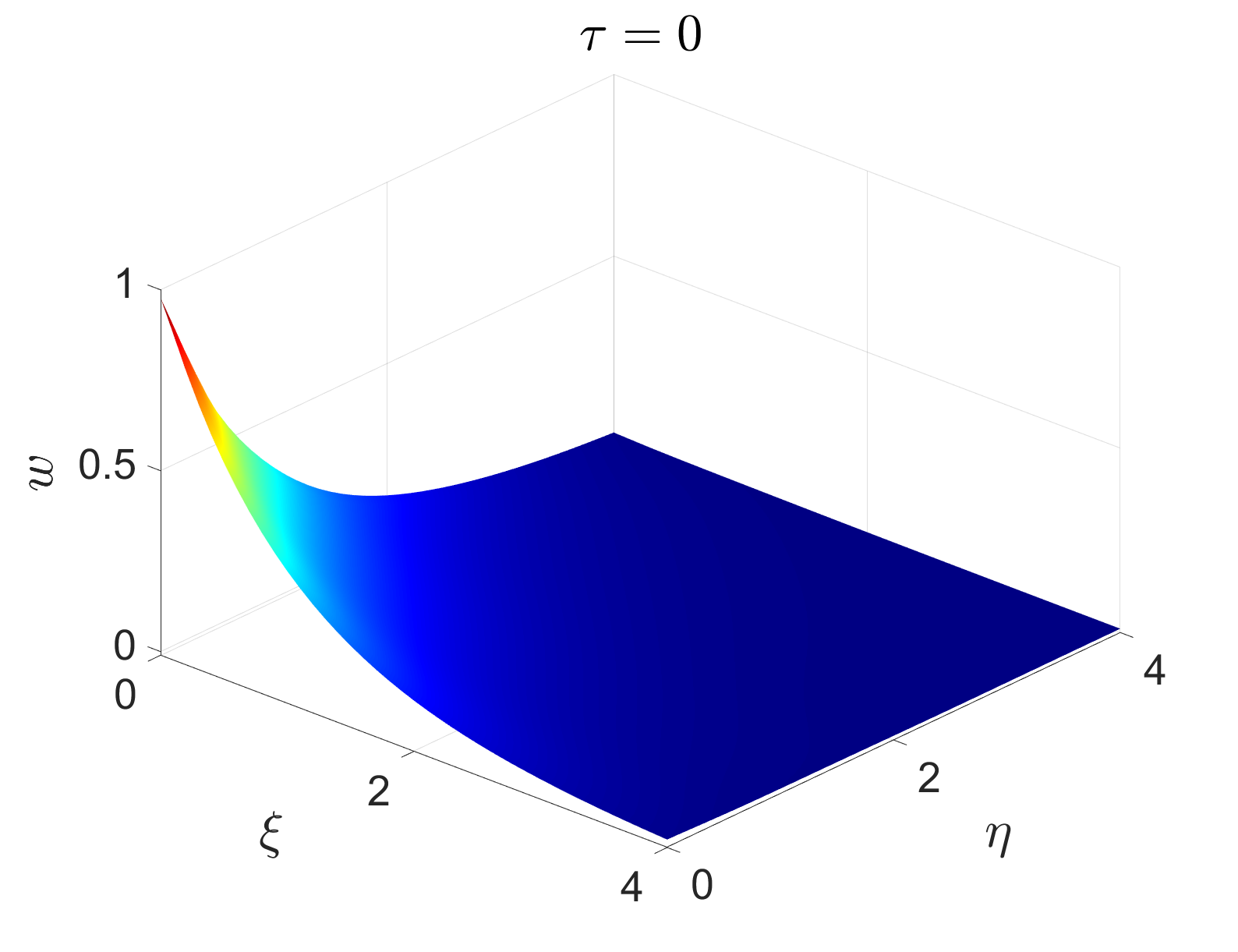}
    \label{fig:diffusion2DMN_a}}
    \subfigure[]{\includegraphics[width=0.45\linewidth]{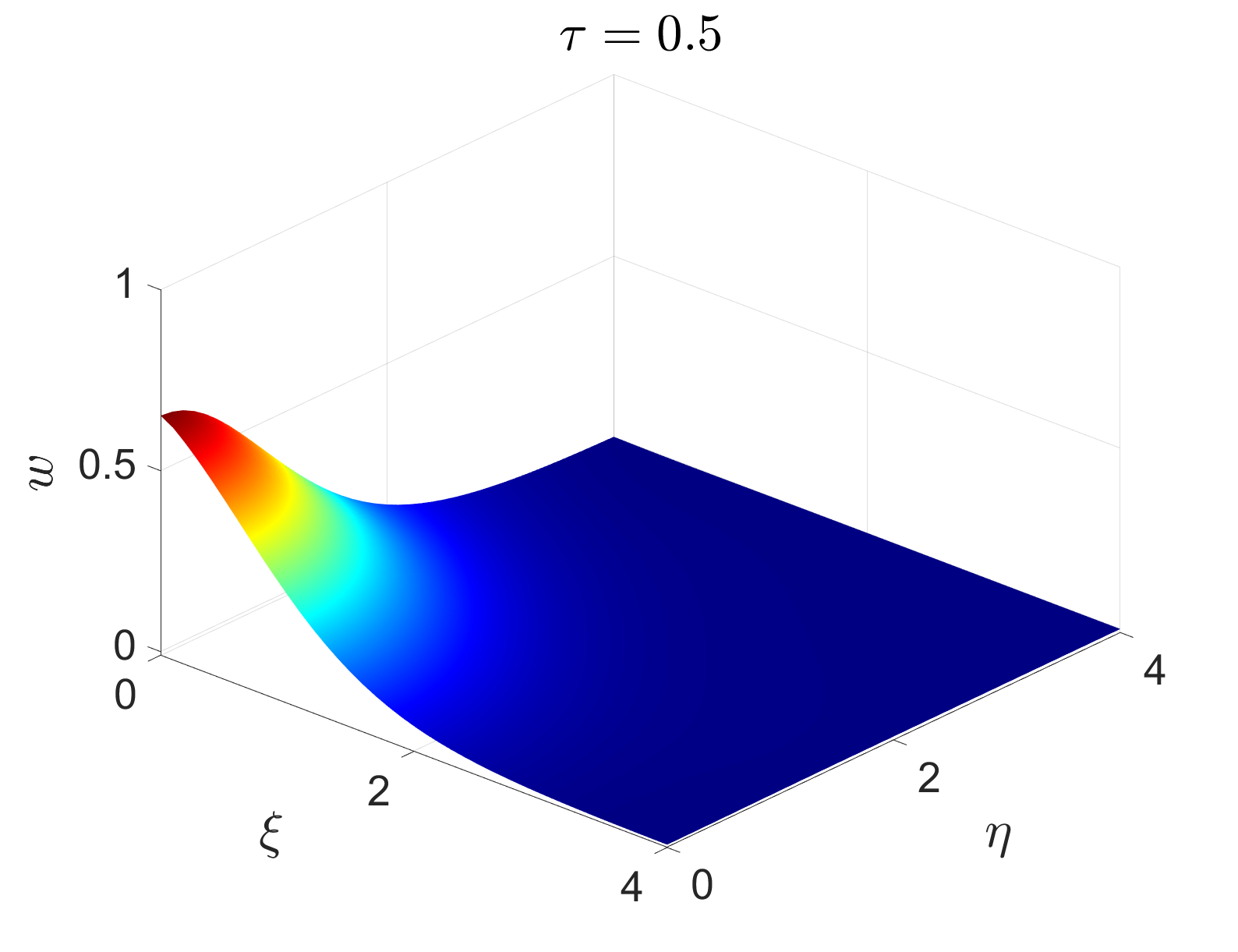}
    \label{fig:diffusion2DMN_b}}
    \subfigure[]{\includegraphics[width=0.45\linewidth]{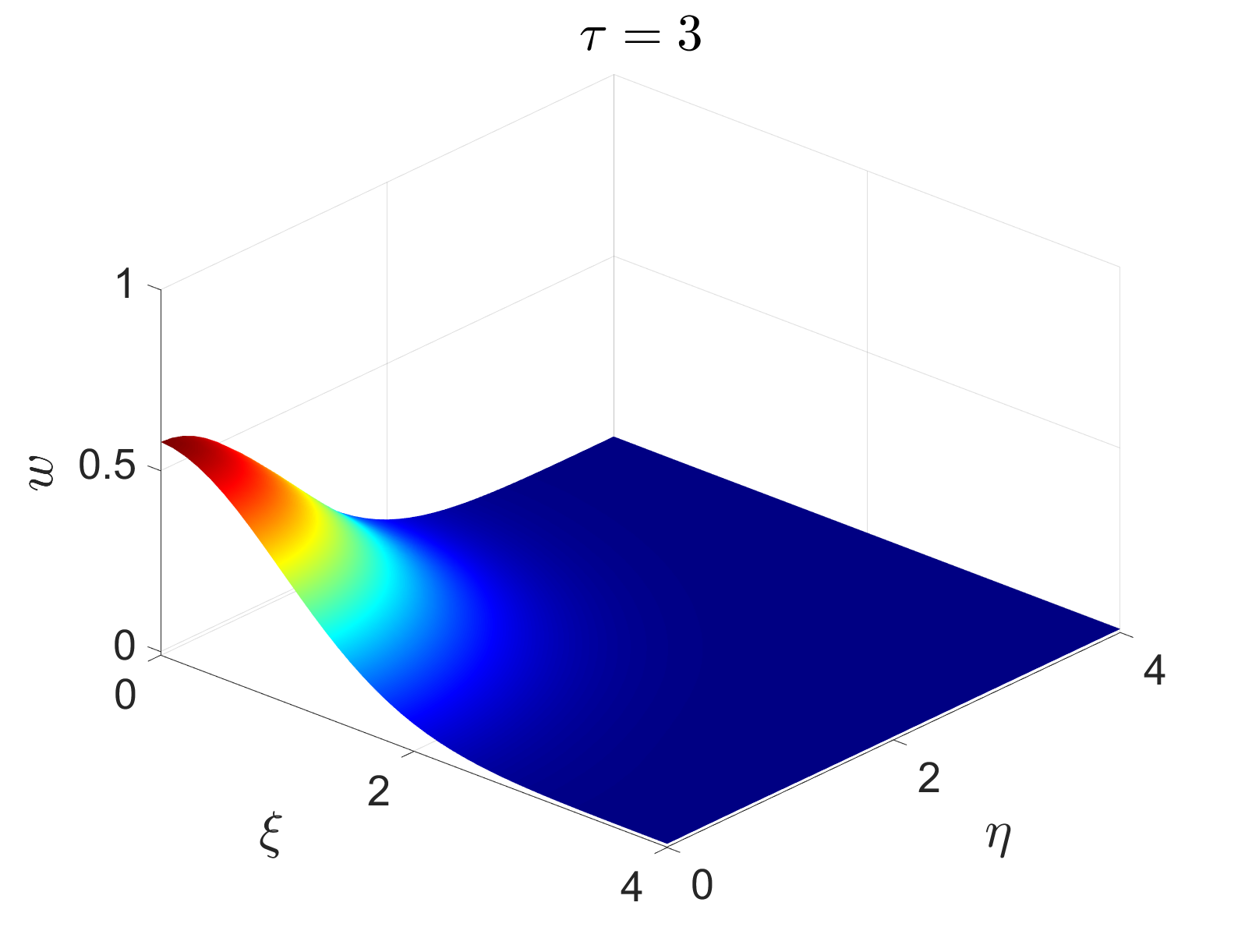}
    \label{fig:diffusion2DMN_c}}
    \subfigure[]{\includegraphics[width=0.45\linewidth]{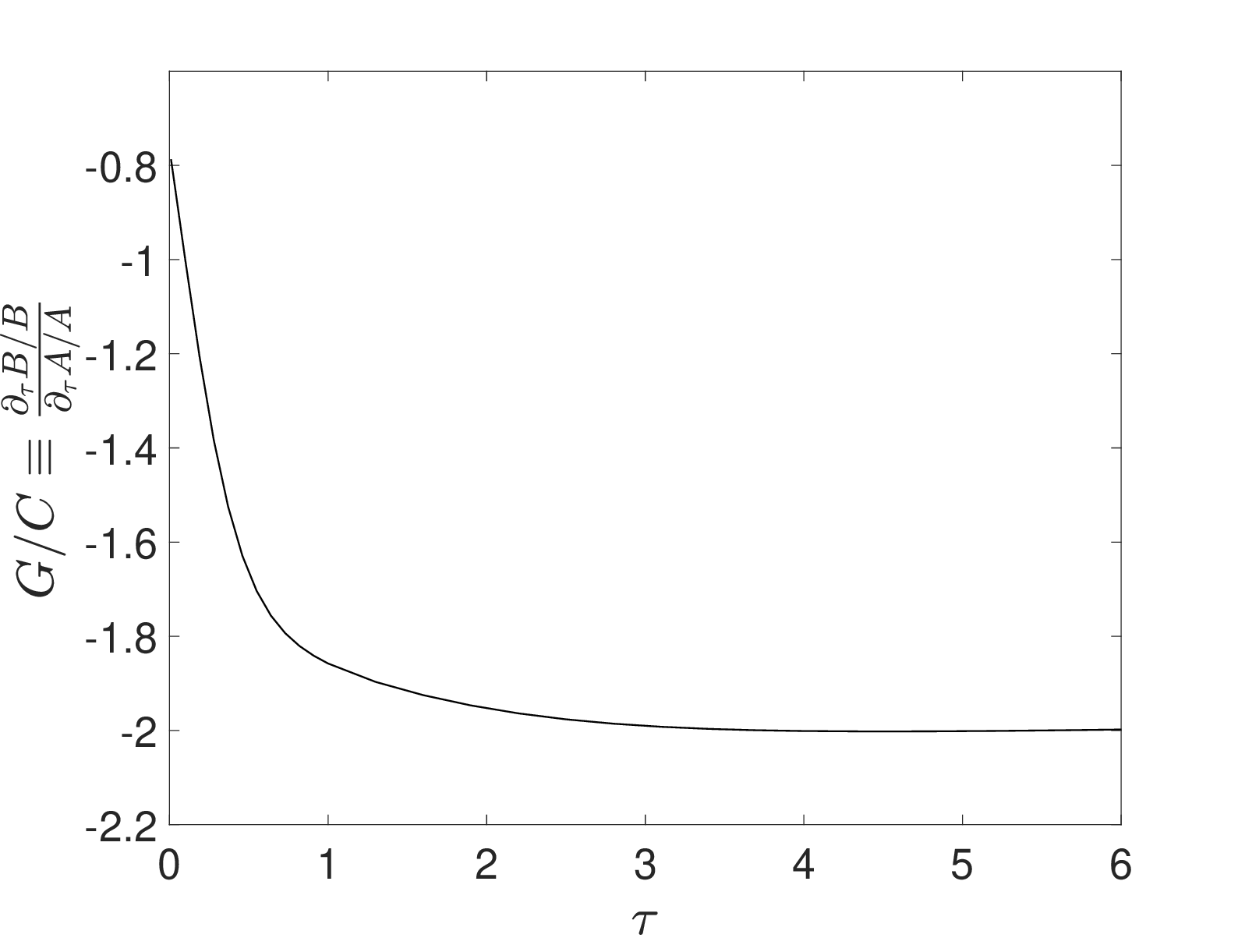}
    \label{fig:diffusion2DMN_d}}
    \caption{PINN prediction of the MN dynamics for the 2D diffusion equation.
    (a)-(c) Snapshots of the rescaled solution $w(\xi,\eta,\tau)$ at rescaled times $\tau=0,0.5,3$.
    The rescaled solution converges to a stationary, self-similar profile for $\tau>3$.
    The neural network $\mathcal{N}_w(\boldsymbol{\xi},\tau)$ is trained using $70,321$ collocation points in the rescaled spatial-temporal $(\xi,\eta)-\tau$ domain.
    (d) Evolution of blow-up rates $G/C$ ratio over rescaled time, $\tau$.
    As $\tau \to \infty$, the ratio $C/G$ converges to a steady-state value, $\lim_{\tau \to \infty} {G/C} = -2$, which is used to derive the self-similar exponents of the diffusion equation.
    }
    \label{fig:diffusion2DMN}
\end{figure}

\noindent the right-hand side operator $\mathcal{D}_{\bf{x}} u \equiv \partial_{xx}^2 u + \partial_{yy}^2 u $ satisfies the scaling relation: 

\begin{equation} \label{eq:scalingdif2d}
    \mathcal{D}_{\bf{x}} \left( B w \left( \frac{x}{A}, \frac{y}{A} \right) \right)  = BA^{-2} \mathcal{D}_{\boldsymbol{\xi}}w, \quad \textrm{ with } \boldsymbol{\xi} \equiv (\xi, \eta) = \left( \frac{x}{A}, \frac{y}{A} \right).
\end{equation}

Again, introducing the ansatz: 

\begin{equation}
    u(x,t)= B(\tau) w \left( \frac{x}{A(\tau)}, \frac{y}{A(\tau)}, 
    \tau(t)\right),
\end{equation}

\noindent and substituting in Eq.~(\ref{eq:diffusion2d}) using Eq.~(\ref{eq:scalingdif2d}) boils down to: 

\begin{equation} \label{eq:diffusion2drescaled}
    \partial_\tau w + \frac{\partial_\tau B}{B} w -\frac{\partial_\tau A}{A} \left( \xi \partial_\xi w + \eta \partial_\eta w \right) = \partial_{\xi \xi}^2 w + \partial_{\eta \eta}^2 w.
\end{equation}
Notice that here, given the isotropy of the spatial
operator, we have assumed that the two spatial
directions scale in the same way.
We determine the blow-up rates, $G(\tau) \equiv \frac{\partial_\tau B}{B}$ and $C(\tau) \equiv \frac{\partial_\tau A}{A}$ by imposing two algebraic template condititions, adopting the rationale applied for the 1D diffusion equation:

\begin{equation} \label{eq:templatedif2D}
    \int_{\eta_{min}}^{\eta_{max}} \int_{\xi_{min}}^{\xi_{max}} \left( w - T \right) T \textrm{d} \xi \textrm{d} \eta = 0 \quad \textrm{ and } \quad \int_{\eta_{min}}^{\eta_{max}} \int_{\xi_{min}}^{\xi_{max}} \left( w - T \right) \left( \xi \partial_\xi T + \eta \partial_\eta T\right) \textrm{d} \xi \textrm{d} \eta = 0, 
\end{equation}

\noindent with $T(\xi, \eta)$ a template function.
For our illustrative example, we solve Eq.~(\ref{eq:diffusion2drescaled}) with initial condition: 

\begin{equation}
\label{eq:initial_diffusion2d}
    w(\xi,\eta, 0) = \exp \left(- (| \xi| + |\eta|) \right)
\end{equation}

\noindent in the domain $(\xi,\eta) \in (0,4) \times (0, 4) $ with symmetry boundary conditions at $\xi=0$, $\eta=0$, and Dirichlet boundary conditions, $w(\xi, \eta, \tau>0)=0$ at boundaries $\xi=4$ and $\eta=4$. 
The template condition $T(\xi,\eta)$ is selected as the initial condition $w(\xi,\eta,0)$.
The neural network $\mathcal{N}_w (\boldsymbol{\xi}, \tau)$ consists of three hidden layers with 40 neurons each, and uses hyperbolic activation functions. 
Similarly to the 1D diffusion equation, the auxiliary network $\mathcal{N}_p (\tau)$ for the prediction of blow-up rates, $G(\tau)$ and $C(\tau)$ has a single hidden layer with five neurons. 
%
%
%
The neural network $\mathcal{N}_w(\boldsymbol{\xi},\tau)$ is trained using $70,321$ collocation points in the rescaled spatial-temporal $(\xi,\eta)-\tau$ domain.
Figures~\ref{fig:diffusion2DMN_a}-\ref{fig:diffusion2DMN_c} depict snapshots of the rescaled solution $w$ at rescaled $\tau=0,0.5,3$. 
The neural network predicted that the rescaled solution $w$ practically converges to a stationary self-similar profile after $\tau \approx 3$. 
In addition, Figure~\ref{fig:diffusion2DMN_d} shows the evolution of $\mathcal{N}_p(\tau)$ predicted ratio of blow-up rates $G/C \equiv \frac{\partial_\tau B /B}{\partial_\tau A /A}$ also shows convergence to a stationary value of approximately $-2$.
In the same way as before this leads to a scaling
of the width with $A \propto (t+t_0)^{1/2}$, while 
 $B \propto (t+t_0)^{-G_{steady}/(2 C_{steady})}$
 which leads to an amplitude exponent of $\approx -1$.
%

\subsection{Case study 3: Second kind self-similar solutions of the axisymmetric 2D Porous Medium Equation}
We now examine the application of MN dynamics to identify self-similar solutions and self-similar exponents of the axisymmetric 2D porous medium equation: 

\begin{equation} \label{eq:pme2d}
    \partial_t u = \partial^2_{rr}u^2 + \frac{1}{r}\partial_r u^2.
\end{equation}

In this problem, the self-similar exponents cannot be determined {\it{a priori}} using dimensional analysis; this corresponds to a case of {\it self-similarity of the second kind}.

\begin{figure}[ht!]
    \centering
    \subfigure[]{\includegraphics[width=0.45\linewidth]{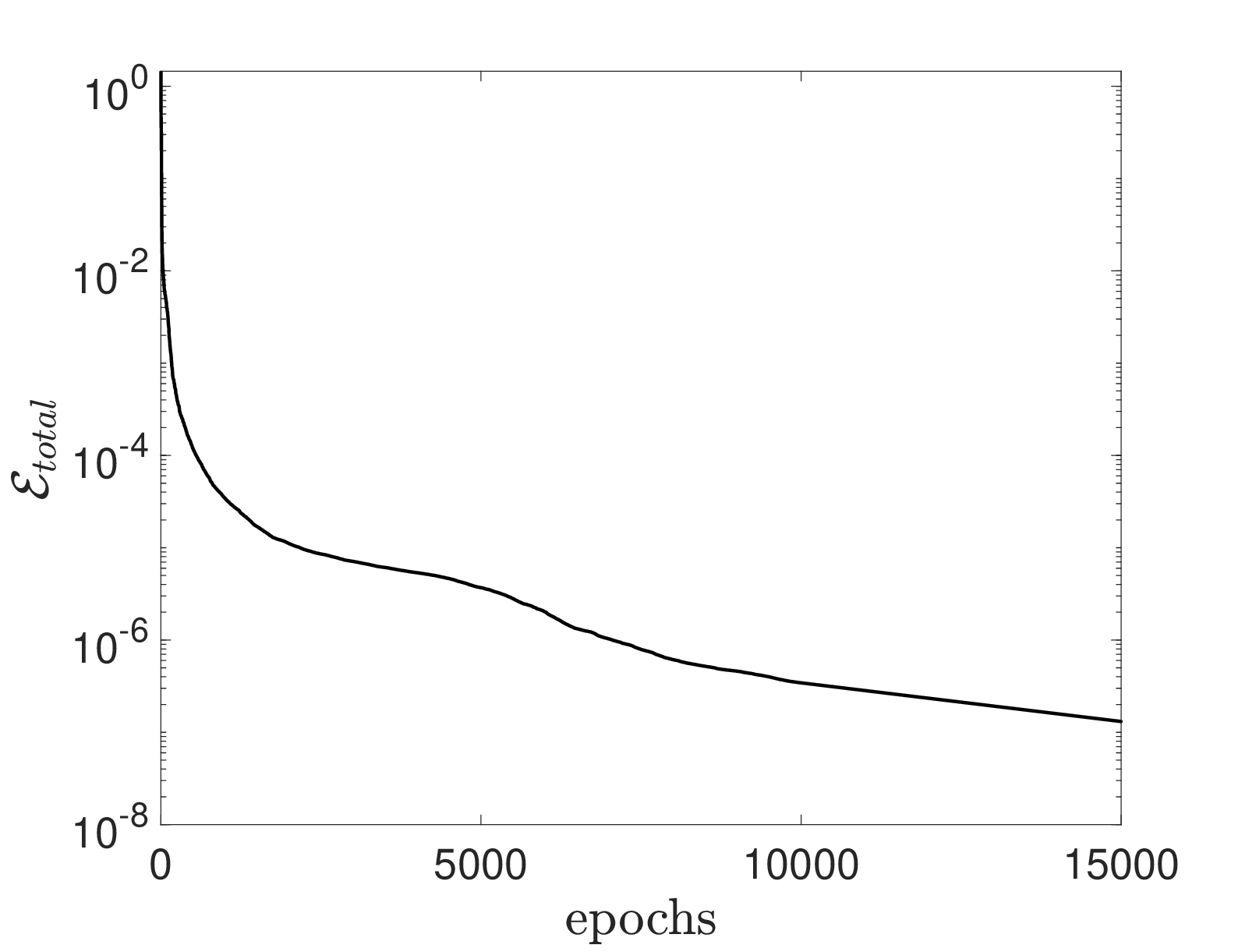}
    \label{fig:pme2DMN_a}}
    \subfigure[]{
    \includegraphics[width=0.45\linewidth]{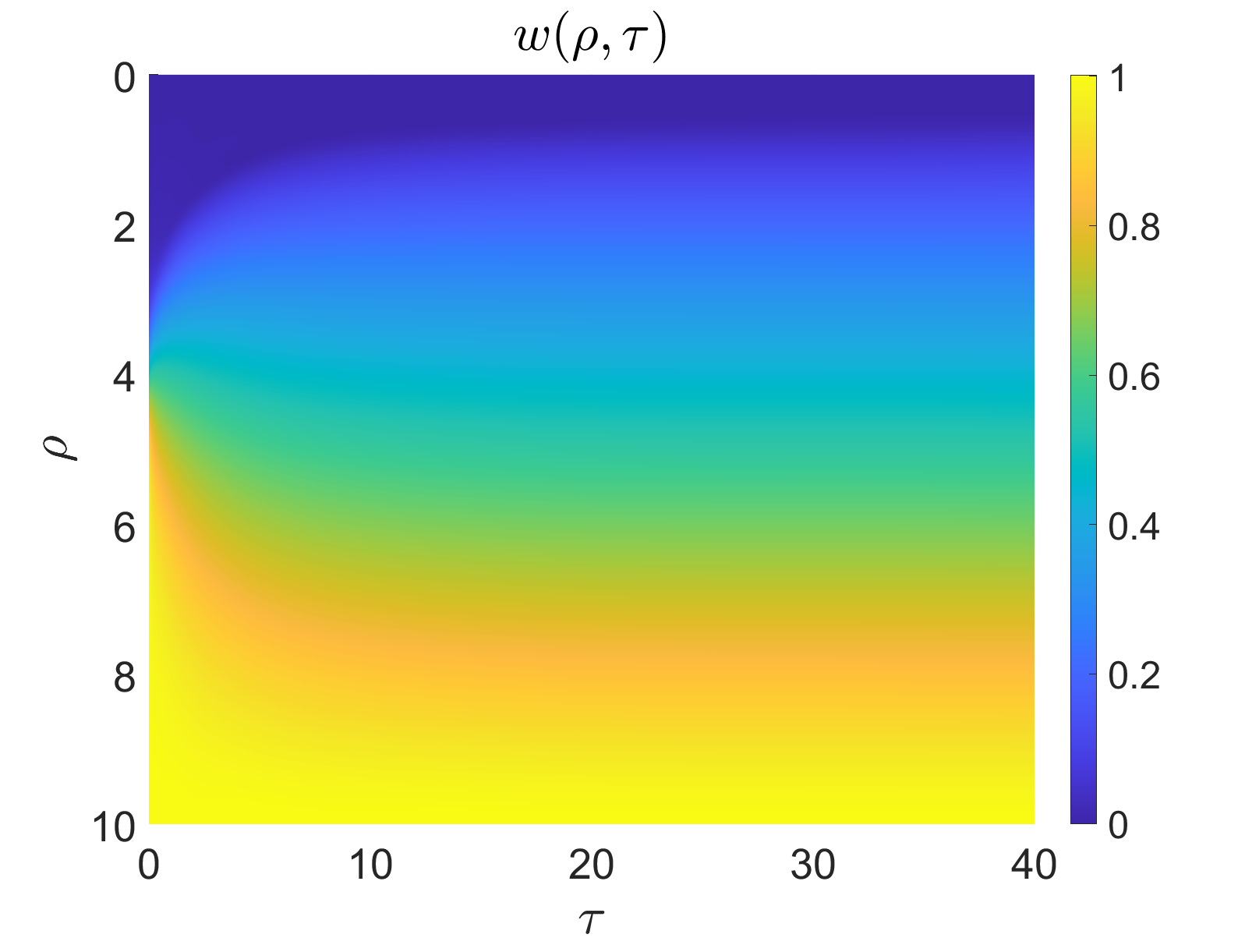}
    \label{fig:pme2DMN_b}}
    \subfigure[]{\includegraphics[width=0.45\linewidth]{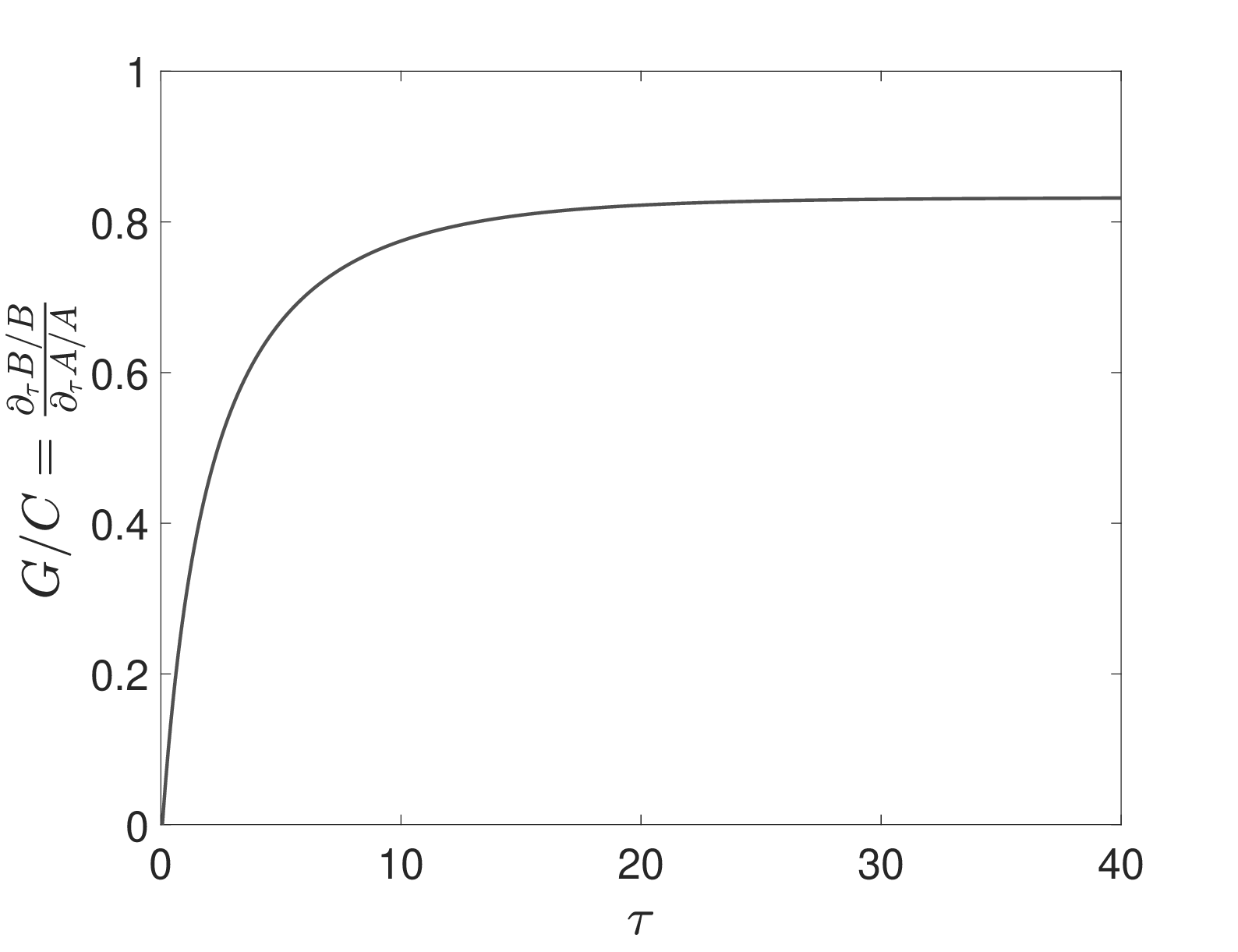}
    \label{fig:pme2DMN_c}}
    \subfigure[]{
    \includegraphics[width=0.45\linewidth]{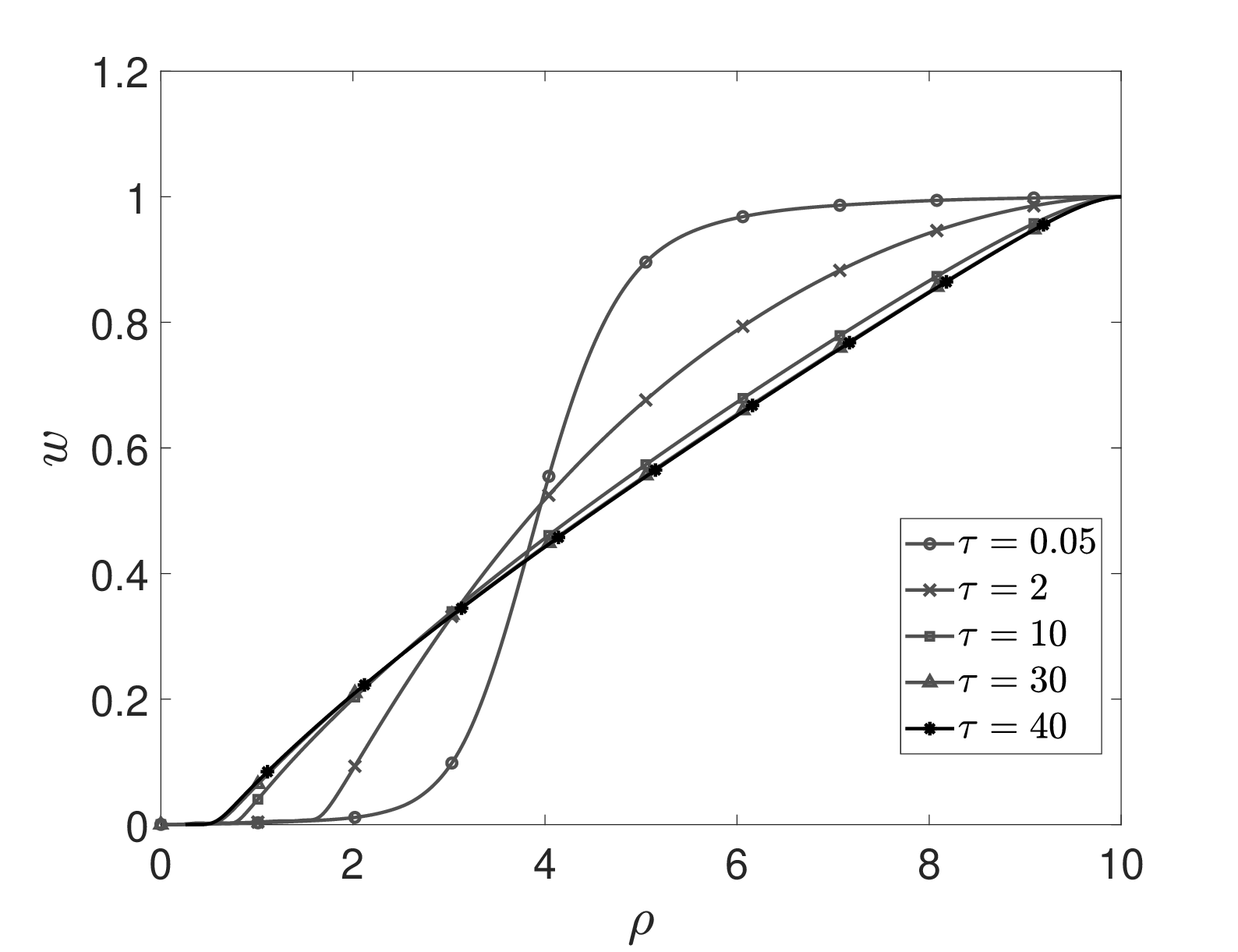}
    \label{fig:pme2DMN_d}}
    \caption{PINN prediction of the MN dynamics for the axisymmetric 2D porous medium equation.
    (a) Evolution of total loss $\mathcal{E}_{loss}$ during training for the rescaled axisymmetric 2D porous medium equation.
    (b) Spatio-temporal evolution of the rescaled solution $w(\rho,\tau)$. 
    The neural network $\mathcal{N}_w(\rho,\tau)$ is trained using 49,900 collocation points in the $\rho-\tau$ domain.
    (c) Evolution of the blow-up rate ratio $G/C$ over rescaled time $\tau$.
    As $\tau \to \infty$, the ratio converges to a steady-state value: $\lim_{\tau \to \infty} {G/C} = 0.831677$, which is used to compute the self-similar exponents of the axisymmetric 2D porous medium equation.
    (d) Snapshots of the rescaled solution, $w(\rho,\tau)$, at various (rescaled) time instances.
    Convergence to the final self-similar profile is effectively achieved by $\tau \approx 25$. 
    \label{fig:pme2DMN}
    }
\end{figure}

\noindent The nonlinear operator $\mathcal{D}_r u \equiv \partial^2_{rr}u^2 + \frac{1}{r}\partial_r u^2 $ satisfies the scaling relation:

\begin{equation} \label{eq:scalingpme}
    \mathcal{D}_r \left( B w \left(\frac{r}{A}\right) \right) = B^2 A^{-2} \mathcal{D}_\rho w, \quad \textrm{ with } \rho = \frac{r}{A}.
\end{equation}

Using the scaling ansatz: 

\begin{equation}
    u(r,t)=B(\tau) w \left( \frac{r}{A(\tau)}, \tau(t) \right)
\end{equation}

\noindent and substituting into Eq.~(\ref{eq:pme2d}) with the scaling relation (\ref{eq:scalingpme}), we obtain the governing equation in the renormalized coordinates:

\begin{equation} \label{eq:pmeMN}
    \partial_\tau w - \frac{\partial_\tau B}{B}w + \frac{\partial_\tau A}{A} \rho \partial_\rho w = \partial^2_{\rho \rho} w^2 + \frac{1}{\rho} \partial_\rho w^2, \quad \textrm{ with } \partial_t \tau = BA^{-2}.
\end{equation}

The unknown blow-up rates $G(\tau) \equiv \partial_\tau B/B$ and $C(\tau) \equiv \partial_\tau A/A$ are computed through two additional algebraic (template) conditions. 
We employ the template conditions proposed in \cite{aronson2001going}. 
The first is an orthogonality condition between the solution, $w$, and a smooth step-like template function, $T(\rho)$:
\begin{equation} \label{eq:templateorthog}
    \int_{\rho_{min}}^{\rho_{max}} w T \textrm{d} \rho =0,
\end{equation}
\noindent where:
\begin{equation}
 T(\rho) = -1 + \frac{2}{1+\exp\left(-2.5(\rho-7)\right)}.   
\end{equation}


%
The second constraint is a pinning condition that fixes the amplitude of the solution at a specific spatial location $p$:
\begin{equation} \label{eq:templatepinning}
    w(p,\tau)=1.
\end{equation}

\noindent The computational domain is $\rho \in [0,10]$.
We impose a zero-flux boundary condition at $\rho=10$ and a Dirichlet boundary condition $w(0,\tau)=0$.
The initial condition for $w$ is given by:

\begin{equation} \label{eq:initial_pme}
w(\rho,0) = \frac{1}{1+\exp\left(-2.5(\rho-4)\right)}
\end{equation}

\noindent We solve the resulting system of DAEs using a PINN approach.
In particular, the network $\mathcal{N}_w(\rho,\tau)$ consists of two hidden layers with 40 neurons each and employs the hyperbolic tangent activation function. 
The blow-up rates $G(\tau)=\partial_\tau B/B$ and $C(\tau)=\partial_\tau A/A$ are predicted by a second neural network, $\mathcal{N}_p(\tau)$, which takes $\tau$ as its only input. 
The neural network $\mathcal{N}_w(\rho,\tau)$ is trained using 49,900 collocation points in the $\rho-\tau$ domain.
The overall loss-function minimized during training includes: $\mathcal{E}_{pde}$ from Eq.~(\ref{eq:pmeMN}), $\mathcal{E}_{bc}$ from boundary conditions, $\mathcal{E}_{ic}$ from the initial condition Eq.~(\ref{eq:initial_pme}), and $\mathcal{E}_{alg}$ from the algebraic template constraints in Eqs.~(\ref{eq:templateorthog}-\ref{eq:templatepinning}).
Training is performed using the L-BFGS optimizer and the convergence of the total loss $\mathcal{E}_{total}$ is shown in Figure~\ref{fig:pme2DMN_a}.
The spatio-temporal solution is depicted in Figure~\ref{fig:pme2DMN_b}.
As demonstrated earlier for the 1D diffusion case, the ratio of these blow-up rates is used to compute the self-similar exponents of the axisymmetric 2D porous medium equation. 

We thus evolve the system up to $\tau=40$ where it practically reaches a steady state and from that
we infer the values of $G_{steady}$ and $C_{steady}$.
These, in turn, lead to the temporal scaling of
$A \propto (t+t_0)^{C_{steady}/(2 G_{steady}-C_{steady})}$,
with the relevant exponent being $\approx 0.85593$,
while
and $B \propto (t+t_0)^{G_{steady}/(2 G_{steady}-C_{steady})}$, with an exponent $\approx 0.711856$; again $t_0$ is an integration
constant.
Once again, this result shows excellent agreement with the theoretical value, e.g., $\alpha=0.85633$ of the exponent 
for the temporal dependence of $A$, 
confirming the accuracy of the proposed PINN-assisted, MN framework for capturing the self-similar solution and exponents in the axisymmetric 2D porous medium equation.

%

\subsection{Case study 4: A problem with scale {\em and} translational invariance - the Burgers equation}

Consider the one-dimensional Burgers equation, which
has also been recently explored in the context of
PINNs (although for self-similarity purposes) 
in~\cite{shahab2025neuralnetworksbifurcationlinear}:
\begin{equation} \label{eq:burgers}
    \partial_tu = \nu \partial^2_{xx} u - u\partial_xu \equiv \mathcal{D}_x(u).
\end{equation}
Here, the value of viscosity is set $\nu=0.025$.
\noindent We apply the scaling ansatz:
\begin{equation}
    u(x,t) = Bw\left( \frac{x-c}{A} \right),
\end{equation}

\begin{figure}[ht!]
    \centering
    \subfigure[]{
    \includegraphics[width=0.45\linewidth]{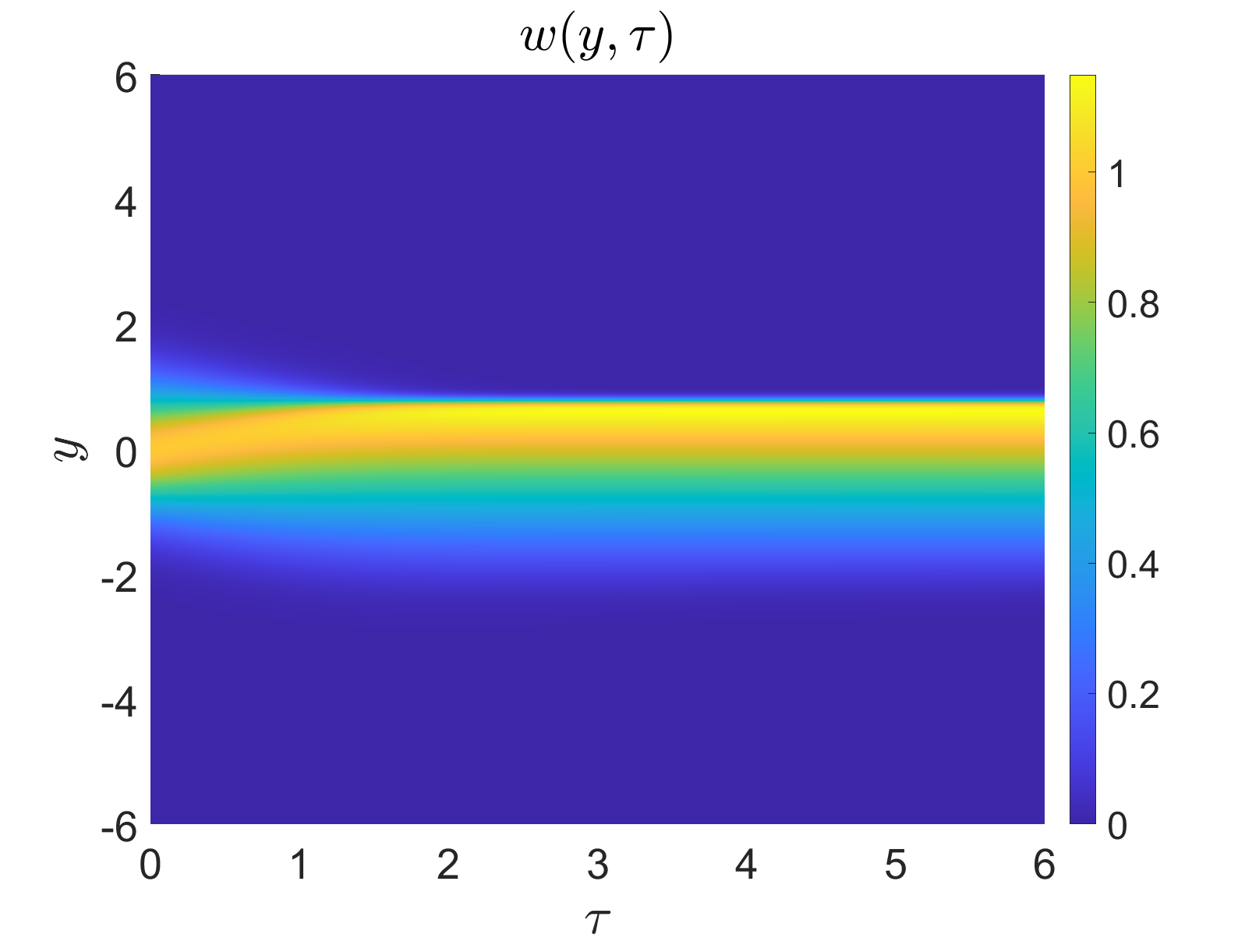}
    \label{fig:burgersMN_a}}
    \subfigure[]{\includegraphics[width=0.45\linewidth]{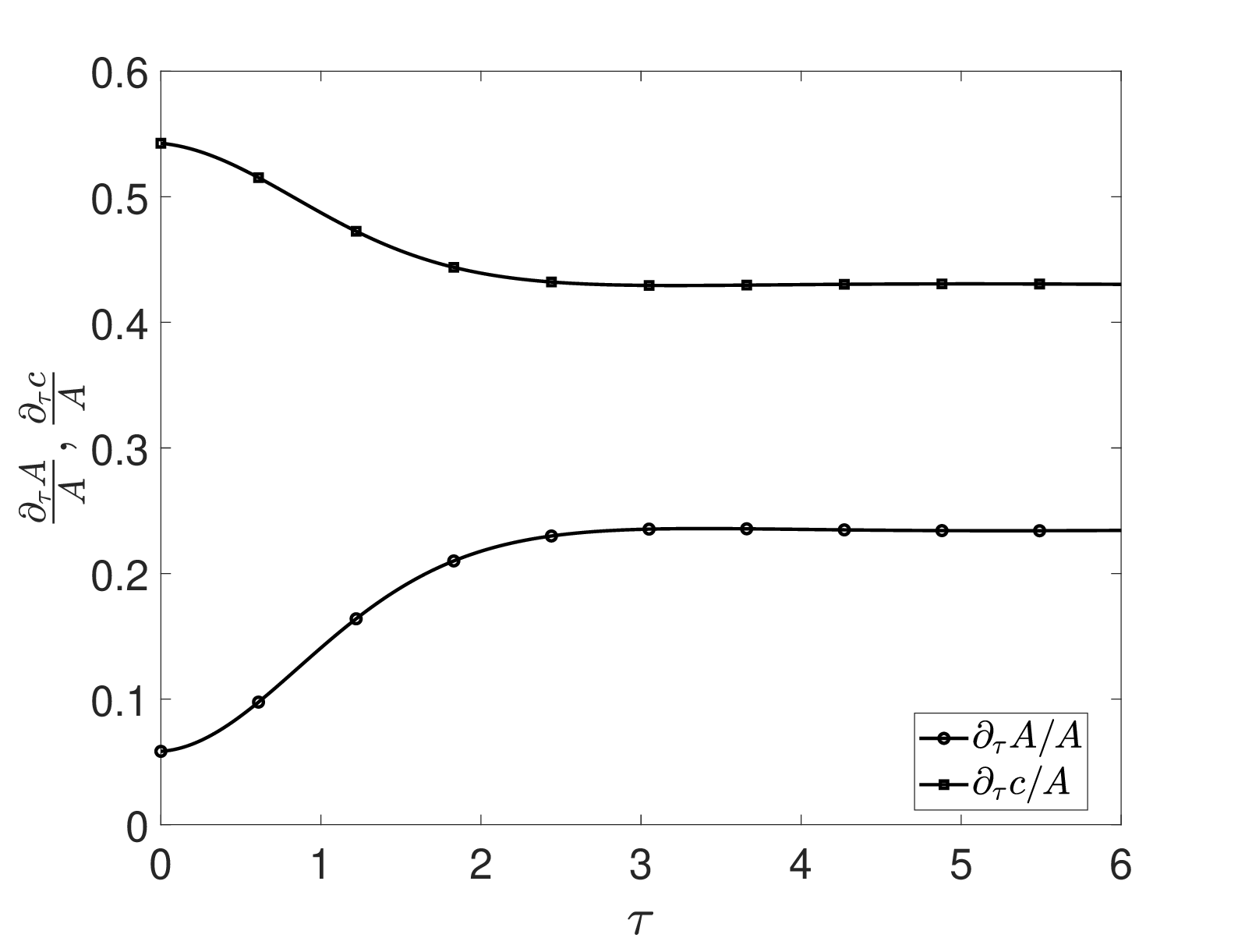}
    \label{fig:burgersMN_b}}
    \subfigure[]{
    \includegraphics[width=0.45\linewidth]{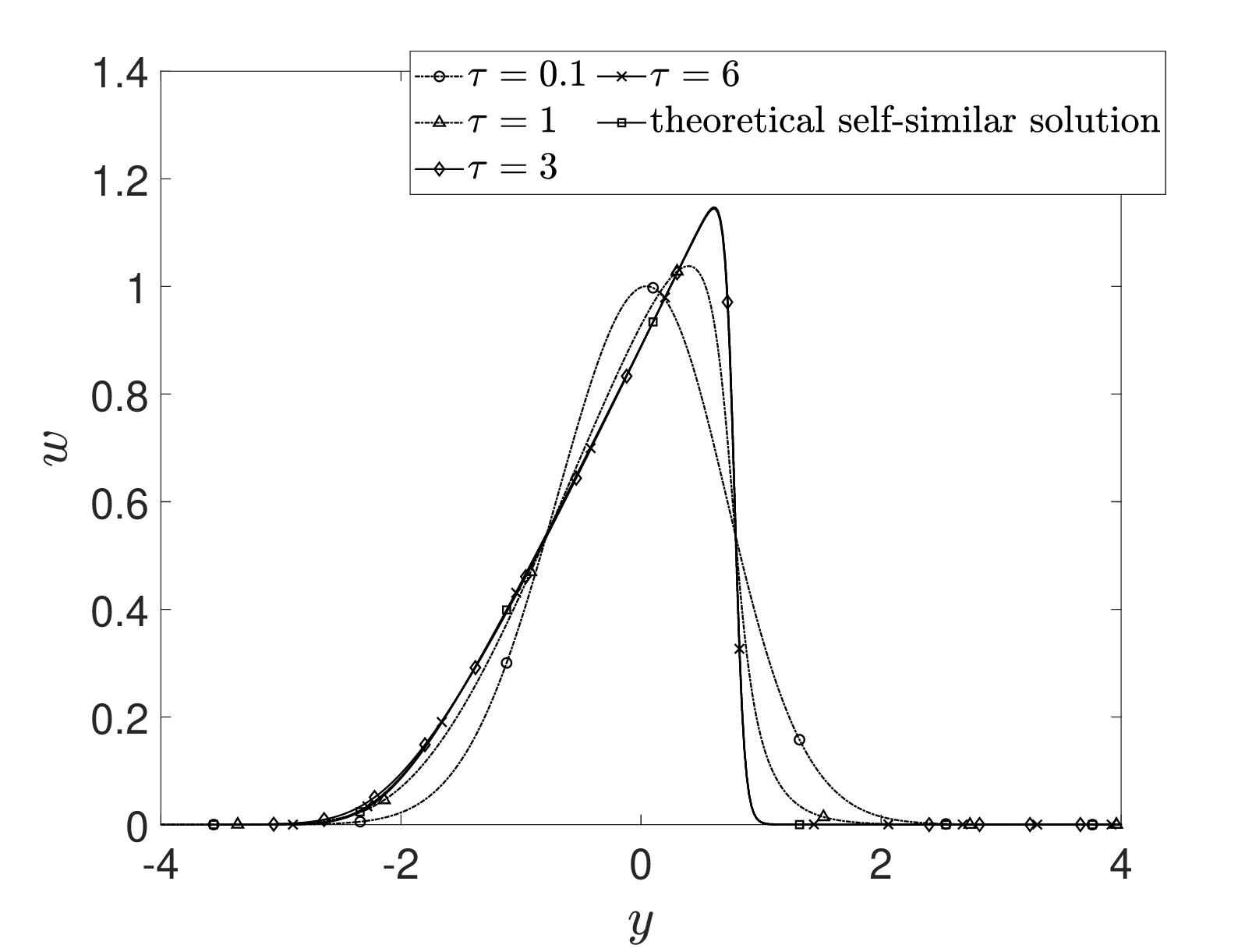}
    \label{fig:burgersMN_c}}
    \subfigure[]{
    \includegraphics[width=0.45\linewidth]{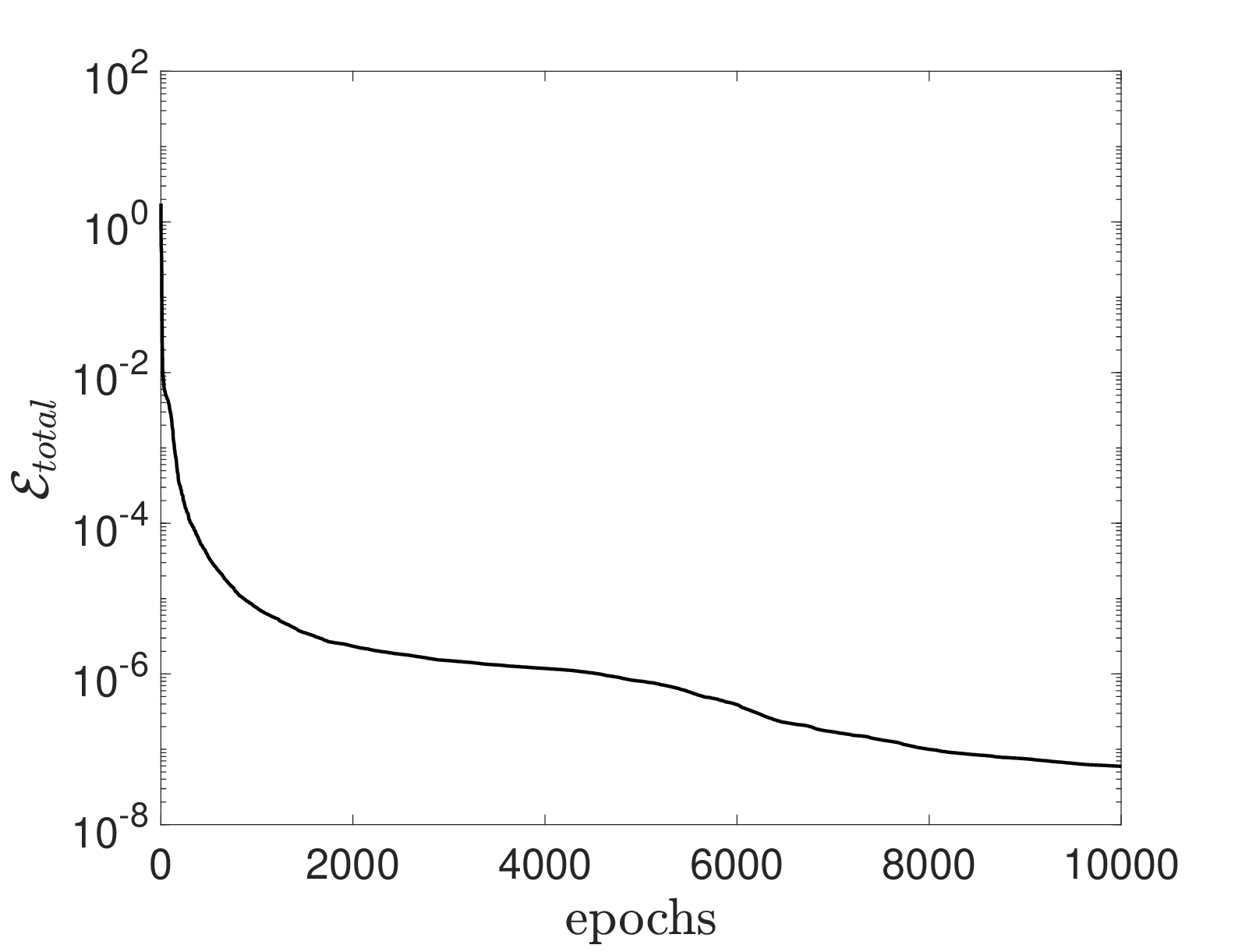}
    \label{fig:burgersMN_d}}
    \caption{PINN prediction of the MN dynamics for the Burgers equation. The value of viscosity $\nu$ in Eq.~(\ref{eq:burgers}) is $\nu=0.025$.
    (a) Convergence of the total loss $\mathcal{E}_{loss}$ during training for the rescaled/co-traveling Burgers equation.
    The neural network $\mathcal{N}_w(y,\tau)$ is trained using 179,700 collocation points in the $y-\tau$ domain.
    (b) Spatio-temporal evolution of the dynamically rescaled solution $w(y,\tau)$. 
    (c) Evolution of $\partial_\tau A/A$ and $\partial_\tau c/A$ over rescaled time $\tau$.
    For long $\tau$ values, ($\tau>3$), the rescaled solution $w$ converges to a stationary, self-similar profile.
    (d) Snapshots of rescaled PINN predicted solution $w$ at different $\tau$ values. 
    The converged self-similar solution $w$ (after practically $\tau \approx 3$) and the theoretical self-similar solution \cite{whitham2011linear} show excellent agreement.
    \label{fig:burgersMN}
    }
    
\end{figure}

\noindent which yields the following relation for the operator $\mathcal{D}_x$:

\begin{equation}
    \mathcal{D}_x \left( B w \left( \frac{x-c}{A} \right) \right) = B A^{-2} \nu \partial^2_{yy}w - B^2A^{-1} w \partial_yw, \quad \textrm{ with } y=\frac{x-c}{A}.
\end{equation}

This means that the operator $\mathcal{D}_x$ is invariant under the transformation when the condition $ AB=1 $ holds.
Under this scaling, we obtain:

\begin{equation}
 \mathcal{D}_x \left(\frac{1}{A} w \left( \frac{x-c}{A} \right) \right) = A^{-3} \mathcal{D}_y(w).   
\end{equation}

Now, substituting the dynamic scaling form: 
\begin{equation}
    u(x,t)=\frac{1}{A(\tau)} w \left( \frac{x-c(\tau)}{A(\tau)}, \tau(t) \right)
\end{equation}

\noindent into Eq.~(\ref{eq:burgers}) leads to the dynamically rescaled equation:

\begin{equation} \label{eq:burgersMN}
    \partial_\tau w = \nu \partial^2_{yy} w - w \partial_yw + \frac{\partial_\tau A}{A} \left( w + y \partial_y w \right) + \frac{\partial_\tau c}{A} \partial_y w, \quad \textrm{ with } y=\frac{x-c}{A}.
\end{equation}

\noindent The two $\tau$-dependent unknown quantities, $\partial_\tau A/A$ and $\partial_\tau c/A$ are determined by imposing two template constraints. 
These arise by minimizing:
\begin{equation}
 E \equiv \int_{y_{min}}^{y_{max}} \left( w - \frac{1}{A} T \left( \frac{y-c}{A} \right) \right)^2 \textrm{d} y,   
\end{equation}

evaluated at $A=1$ and $c=0$. $T(y)$ is a chosen template function.
This leads to the following constraints: 

\begin{equation} \label{eq:burgerstemplate}
    \int_{y_{min}}^{y_{max}} \left(w-T\right) \left( T + y \partial_y T \right) \textrm{d} y =0 \quad \textrm{ and } \quad \int_{y_{min}}^{y_{max}} \left(w-T\right) \partial_y T  \textrm{d} y = 0
\end{equation}

Here, we select $T(y)=\exp\left(-y^2\right)$.
Equation~(\ref{eq:burgersMN}) is solved in the domain $y \in [-6,6]$ with zero-flux boundary conditions at both ends. 
The initial condition is $w(y,0)=\exp\left(-y^2\right)$.

For this implementation, the neural network $\mathcal{N}_w(y,\tau)$ consists of three hidden layers with 20 neurons each, using the hyperbolic tangent activation function. 
The auxiliary network $\mathcal{N}_p(\tau)$ which predicts $\partial_\tau A/A$ and $\partial_\tau c/A$ contains one hidden layer with 4 neurons.
The neural network $\mathcal{N}_w(y,\tau)$ is trained using 179,700 collocation points in the $y-\tau$ domain.
Training is performed using the L-BFGS optimizer and the convergence of the total loss $\mathcal{E}_{total}$ is shown in Figure~\ref{fig:burgersMN_a}.
Figure~\ref{fig:burgersMN_b} illustrates the spatio-temporal evolution of the dynamically rescaled solution $w(y,\tau)$, which converges to a stationary self-similar profile as long $\tau$ values (here $\tau=6$).
While the initial data steepens into a sharp front, the front neither propagates nor increases in amplitude - unlike the solution of the original Burgers equation in the $x-t$ domain. Indeed, this demonstrates that we have converged
to a stationary (indeed, stable, since it was
reached through evolution dynamics) self-similar 
traveling solution. 
To evaluate the accuracy of our PINN-assisted framework, we compare the rescaled solution $w(y,\tau)$ at $\tau=6$ - where the solution has effectively converged- with the analytical self-similar solution of the Burgers equation \cite{whitham2011linear}:

\begin{equation} \label{eq:burgersanalytic}
    u(x,t_0)=\sqrt{ \frac{\nu}{\pi t^*}} \frac{\left[ \exp \left(A^*/\left( 2 \nu \right) \right) -1 \right] \exp \left[ - \left(x - c^* \right)^2 / \left(4 \nu t^* \right) \right]}{1+ \left[ \exp \left(A^*/\left(2 \nu \right) \right)-1 \right]/2 \cdot \textrm{erfc} \left( \left(x-c^* \right)/\sqrt{4 \nu t^*} \right)}
\end{equation}

The evolution of $\partial_\tau A/A$ and $\partial_\tau c/A$ over rescaled time $\tau$ is shown in Figure~\ref{fig:burgersMN_c}.
As shown in Figure~\ref{fig:burgersMN_d}, there is excellent agreement between the  theoretical self-similar solution  and the PINN-predicted profile $w$, with optimal parameters: amplitude, $A^* \approx 1.7713$, spatial shift: $c^* \approx  -1.8715 $ and time shift: $t^* \approx  2.1666$.


\section{Conclusions and Future Challenges}
\label{sec:conclusions}
In this work we developed a PINN-based framework for solving nonlinear PDEs with continuous symmetries -- specifically translations and scalings -- by reformulating them as high-index DAEs in dynamically rescaled coordinates.
These index-2 DAEs simultaneously encode the evolution of the invariant profile and the symmetry parameters (such as wave speed or
similarity exponents). To solve these systems, we employed Physics-Informed Neural Networks (PINNs), which integrate the PDE residuals and algebraic constraints, arising from symmetry reduction, directly into the loss function. This enables the direct recovery of invariant solutions and their associated group parameters (e.g., wave speed or similarity exponents).

The method was tested on four canonical problems, involving traveling waves and first- and second-kind self-similar behavior. In all cases, PINNs were able to solve the resulting constrained systems, making this a promising approach for studying
scaling and front-propagating phenomena.
One of the key advantages of this symmetry-reduced approach is its ability to handle problems that are otherwise difficult to solve using traditional numerical schemes. Many self-similar and moving-front problems would require adaptive meshes or moving domains to maintain resolution near sharp features. In contrast, the dynamic rescaling used here maintains the structure of steep or concentrating solutions (e.g., solitons or traveling waves) at a constant scale in the transformed coordinates, reducing the need for mesh refinement; notice that this is a potential 
advantage in comparison to PINN methods that 
operate in the original frame within such problems~\cite{shahab2025neuralnetworksbifurcationlinear}.
This also tends to reduce stiffness and removes fast transients (we found that the reduced dynamics often converge to stationary or slowly-varying states, facilitating training) potentially improving numerical stability. 
Additionally, the method in its implementation herein
bypasses the need to derive scaling conditions
beforehand, as well as the technical difficulty
of large equation residues of PINNs at particular
points (when applied
witout pinning/template conditions)~\cite{wang2023asymptotic}. While the latter
leads to extra (potentially custom made for different 
solutions) loss terms that eventually enable the
identification of scaling exponents, these types of
considerations which appeared to be necessary in earlier
implementations, are systematically accounted for by
the formalism
and therefore absent in the present formulation.
Nevertheless, a more  systematic analysis of the conditioning of the transformed problem across various PDEs remains to be done.

Traditional solvers, like finite differences with method-of-lines are typically faster, but become difficult to apply in high-index DAE formulations. They require index reduction or specialized time-integrators, which are not standard in environments like \texttt{MATLAB}, and often require external specialized software. Thus,
the PINN offers a more flexible alternative, especially when dealing with algebraic constraints and dynamic symmetries in complex geometries or under limited software availability.
In this context, the flexibility of PINNs makes them appealing for such applications, not necessarily due to performance, but due to their unified treatment of residuals and constraints and ease of implementation.
Moreover, one of the real strengths of PINNs that
is especially relevant to leverage in future
work is its ability to handle efficiently problems
in higher dimensions, given their optimization 
(rather than fully well-posed solution nature 
as in traditional methods). Relevant comparisons
have recently appeared, e.g., in~\cite{carola}. 

The approach is not without drawbacks. PINNs are known to be sensitive to hyperparameters and network architecture, and optimizer choice plays a crucial role.
We used a simple feedforward network with tanh activations; other architectures may perform better or worse depending on the specific problem, but choosing optimal configurations is not trivial.
We also observed that standard gradient-based optimizers such as ADAM often fail to achieve sufficient accuracy, particularly when the loss includes algebraic constraints that must be satisfied to high precision. The L-BFGS optimizer, while more effective in this context, is expensive. Its memory and computational cost scales poorly with network size, limiting the method to small- and medium-sized networks in practice.
To mitigate convergence issues, we employed warm-up strategies, first fitting initial and boundary conditions before enforcing the full physics and algebraic constraints. We also tested greedy, continuation-inspired training schemes \cite{alvarez2023discrete,alvarez2024nonlinear}, where the solution is first learned on a small subdomain (e.g., close to initial conditions or pinning conditions) and then gradually extended. This corresponds to a zero-order continuation approach, useful in problems where steep gradients or complex transients emerge. Yet in the transformed coordinates, the dynamics were already regular enough that such strategies showed no clear benefit.
%
We acknowledge that the PINN approach,
as implemented here, is not optimally computationally efficient. Training involves minimizing a non-convex loss function over a large parameter space, using second-order optimizers such as L-BFGS. While effective, this approach can be computationally demanding.
Additionally, in the present work, we have sought
to allow the integrator to converge to a {\it stable}
self-similar or traveling (or both) solution.
However, it seems quite relevant
to seek to leverage fixed point
methods (in the self-similar or co-traveling frame)
to identify such solutions. Moreover, the latter
aspect may enable both the more accurate identification
of solutions (in a controllable way) and their use
for the consideration of stability features~\cite{CHAPMAN2022133396}, as well
as the potential identification of unstable,
e.g., multi-hump self-similar solutions, as was the
case, e.g., in~\cite{BUDD1999756}

Looking forward, more efficient alternatives such as parsimonious Random Projection Neural Networks (RPNN), as in \cite{fabiani2023parsimonious}, might offer significant computational advantages, by employing randomized features and linear-algebraic solvers (e.g., SVD, COD) rather than gradient-based optimization.  In particular, RPNN-based PINNs were shown in \cite{fabiani2023parsimonious} to reach accuracies comparable to classical index-1 DAE solvers like \texttt{ode15s} and \texttt{ode23t} in \texttt{MATLAB}.
We plan to investigate RPNN-based PINNs combined with method-of-lines discretization to enable a more direct comparison with classical numerical solvers, while retaining the flexibility and implementation advantages of the PINN framework. This would allow a clearer assessment of performance and accuracy across both paradigms.
Moreover, it would be of particular interest to 
explore models with more complex dynamics
and potential bifurcations of self-similar solutions,
such as e.g., the nonlinear Schr{\"o}dinger~\cite{SulemSulem1999,Fibich2015} and
the Korteweg-de Vries equations~\cite{eva,Chapman_2024}.
Identifying self-similar solutions as steady-states through the use of Newton iteration and pseudo-arclength continuation allows for the computation of both stable, and unstable families of such solutions \cite{siettos2003focusing,chapman2024self,chapman2025multi}.
In summary, this method provides a flexible way to study nonlinear PDEs with symmetries, particularly when standard techniques are difficult to apply. It opens a path toward scalable neural PDE solvers with structure-aware formulations. Yet, more and more systematic explorations
of the nature of the networks, the choice of hyperparameters,
and the selection of the optimization are clearly needed
to identify how to best perform the relevant tasks
for the models of interest.

\section*{Declarations}

\subsection*{Data and code availability}
The data and code supporting the findings of this study will be available upon publication or upon reasonable request.

\subsection*{Acknowledgments}
G.F. and I.G.K. acknowledge the Department of Energy (DOE) support under Grant No. DE-SC0024162,
as well as partial support by the National Science Foundation under Grants No. CPS2223987 and FDT2436738.
C.S. acknowledges partial support by the PNRR MUR projects PE0000013-Future Artificial Intelligence Research-FAIR \& CN0000013 CN HPC - National Centre for HPC, Big Data and Quantum Computing, Gruppo Nazionale Calcolo Scientifico-Istituto Nazionale di Alta Matematica (GNCS-INdAM).

\subsection*{Competing interest statement}
The authors declare that there are no known competing financial interests or personal relationships that could have appeared to influence the work reported in this manuscript.

\bibliographystyle{abbrv}
\bibliography{AA_references.bib}

\clearpage
\appendix

\end{document}